\documentclass[12pt]{amsart}
\usepackage[utf8]{inputenc}

\font\myfont=cmr12 at 20pt

\title[Daubechies' localization operator on Cantor Type Sets II]{\myfont{\textbf{Daubechies' Time-Frequency Localization Operator on Cantor Type Sets II}}}
\author{Helge Knutsen}
\email{helge.knutsen@ntnu.no}
\date{April 2021}
\address{Department of Mathematical Sciences, Norwegian University of Science and \mbox{Technology,} \newline 7034 Trondheim, Norway}

\usepackage{amsthm}
\theoremstyle{plain} 

\newtheorem{theorem}{Theorem}[section]
\newtheorem{corollary}{Corollary}[section]
\newtheorem{lemma}[theorem]{Lemma}
\newtheorem{proposition}{Proposition}[section]

\theoremstyle{definition} 
\newtheorem*{remark}{Remark}
\newtheorem{example}{Example}[section]

\usepackage{amsmath, amssymb} 
\numberwithin{equation}{section} 

\usepackage{mathtools}

\DeclarePairedDelimiter\floor{\lfloor}{\rfloor} 

\usepackage{enumitem} 

\usepackage{calrsfs} 
\DeclareMathAlphabet{\pazocal}{OMS}{zplm}{m}{n}

\usepackage{bm} 

\usepackage[toc,page]{appendix} 

\usepackage{graphicx}
\graphicspath{ {figures/} }

\usepackage{pst-node}
\usepackage{tikz-cd}

\usepackage{geometry}
\usepackage{afterpage}

\usepackage{verbatim} 

\usepackage{subfigure}

\usepackage[hyphens]{url}

\usepackage{csquotes}

\usepackage[colorlinks,urlcolor=black]{hyperref} 

\begin{document}

\maketitle



\begin{abstract}
We study a version of the fractal uncertainty principle in the joint time-frequency representation. Namely, we consider Daubechies' localization operator projecting onto spherically symmetric $n$-iterate Cantor sets with an arbitrary base $M>1$ and alphabet $\mathcal{A}$. We derive an upper bound asymptote up to a multiplicative constant for the operator norm in terms of the base $M$ and the alphabet size $|\mathcal{A}|$ of the Cantor set. For any fixed base and alphabet size, we show that there are Cantor sets such that the asymptote is optimal. In particular, the asymptote is precise for mid-third Cantor set, which was studied in part I \cite{HelgeKnutsen_Daubechies_2020}. Nonetheless, this does not extend to every Cantor set as we provide examples where the optimal asymptote is not achieved. 
\end{abstract}


\section{Introduction}
There are many versions of the uncertainty principle, that all, in some form, state that a signal cannot be highly localized in time and frequency simultaneously. One recent version, described by Dyatlov in \cite{DyatlovFUP2019}, referred to as the \textit{fractal uncertainty principle} (FUP) (first introduced and developed in \cite{Dyatlov_Zahl_2016}, \cite{Bourgain_Dyatlov_2018}, \cite{Dyatlov_Jin_2018}), states that a signal cannot be concentrated near fractal sets in time and frequency. Fractal sets are here defined broadly as families of time and frequency sets $T(h), \Omega(h)\subseteq[0,1]$, dependent on a continuous parameter $0<h\leq 1$, that are so-called $\delta$-regular with constant $C_{R}\geq 1$ on scales $h$ to $1$ (see Definition 2.2. in \cite{DyatlovFUP2019} for details). The FUP is then formulated for this general class when $h\to 0$. While originally formulated in the separate time-frequency representation, we search for an analogous result in the joint representation as the uncertainty principles should be present regardless of our choice of time-frequency representation.

In the separate representation, Dyatlov considers the following compositions of projections onto such time and frequency sets, $\chi_{\Omega(h)}\pazocal{F}_{h}\chi_{T(h)}$, where $\chi_{E}$ is the characteristic function of a subset $E\subseteq \mathbb{R}$, and $\pazocal{F}_{h}$ is a dilated Fourier transform $\pazocal{F}_{h}f(\omega) = \sqrt{h}^{-1}\pazocal{F}f(\omega h^{-1})$.  For signals $f\in L^2(\mathbb{R})$, the FUP then states (see Theorem 2.12 and 2.13 in \cite{DyatlovFUP2019}) that for a fixed $0<\delta<1$ there exists some non-trivial exponent $\beta>\max\{0,1/2-\delta\}$ such that
\begin{align}
    \| \chi_{\Omega(h)}\pazocal{F}_{h}\chi_{T(h)}\|_{\text{op}}=\pazocal{O}(h^{\beta}) \ \ \text{as } h\to 0.
    \label{FUP_statement}
\end{align}
Although no further estimates for the exponent is provided in \cite{DyatlovFUP2019}, Jin and Zhang have made some progress in a recent paper \cite{Jin_Zhang_FUP_2020}: For families of sets $T(h),\Omega(h)\subseteq[0,1]$ that are $\delta$-regular with constant $C_{R}\geq 1$ on scales $h$ to $1$, they have derived an explicit estimate for $\beta>0$ only dependent on $0<\delta<1$ and $C_{R}$.

In the original statement, we might not immediately recognize \eqref{FUP_statement} as an uncertainty principle since the parameter $h$ is also encoded in the Fourier transform. However, if we disentangle $h$ from the Fourier transform, \eqref{FUP_statement} turns out to be a statement regarding localization on the sets $T(h)/\sqrt{h}$ in time and $\Omega(h)/\sqrt{h}$ in frequency. Depending on our choice of $\delta$, the measures $|T(h)/\sqrt{h}|, |\Omega(h)/\sqrt{h}|$ might, in fact, tend to infinity as $h\to 0$.

One prominent subfamily of fractal sets, is the Cantor sets. This includes the infamous mid-third Cantor set, but instead of subdividing into three pieces and keeping two by each iteration, we generalize and subdivide into $M>1$ pieces labeled $\{0,1,\dots,M-1\}$ and keep a fixed alphabet $\mathcal{A}$ of said pieces. With regard to the Cantor set construction, we could just as well speak of iterations $n\to \infty$ rather than an $h$-neighbourhood that tends to zero. For Cantor sets with base $M$, these two quantities are related by $h\sim M^{-n}$. 

Moreover, if $T(h) =\Omega(h)$ correspond to an $n$-iterate Cantor set $\pazocal{C}_n(M,\mathcal{A})$ based in $[0,1]$, their scaled counterparts $T(h)/\sqrt{h} = \Omega(h)/\sqrt{h}$ correspond to a Cantor set based in $[0,L]$ with $L\sim M^{\frac{n}{2}}$. Or equivalently, for such scaled $n$-iterate Cantor sets, the intervals $\{I_j\}_j$ that make up $n$-iterate all satisfy
\begin{align}
    |I_j| \sim \frac{1}{L} \sim M^{-\frac{n}{2}}.
    \label{FUP_condition}
\end{align}
This condition will be a point of reference in the subsequent discussion. 

In the joint representation, we shall consider Daubechies' localization operator, based on the Short-Time Fourier Transform (STFT) with a Gaussian window, that projects onto spherically symmetric subsets of the time-frequency plane. The non-trivial assumption of radial symmetry is effective as the Daubechies operator has a known eigenbasis, the Hermite functions, and explicit formulas for the eigenvalues. With this powerful tool available, we estimate the operator norm when projecting onto the family of generalized spherically symmetric Cantor sets with base $M$ and alphabet $\mathcal{A}$ defined in a disk of radius $R>0$. Other treatments of localization on sparse sets utilizing the STFT can be found in \cite{Fernandez_Galbis_2010} and \cite{Abreu_Speckbacher_2018}, where sparsity is described in terms of "thin at infinity" and the Nyquist density, respectively. 

The remainder of the paper is organized as follows: In section \ref{section_preliminaries}, we describe Daubechies' operator and the Cantor set construction in more detail. New results are divided into two sections, \ref{localization_on_gen_spherically_symmetric_Cantor_set} and \ref{section_indexed_cantor_set}: In the first section, we keep the base $M>1$ and alphabet $\mathcal{A}$ \textit{fixed} throughout the iterations, while in the second, we introduce the notion of an indexed Cantor set where the base and alphabet may \textit{vary}. 

In section \ref{localization_on_gen_spherically_symmetric_Cantor_set}, we let the radius $R = R(n)$ be dependent on the iterates $n$. For base $M$, we initially consider the general case when
\begin{align}
    R(n)\to \infty \ \ \text{as } n\to \infty \ \ \text{while} \ \ \pi R^2(n)\leq M^{n} \ \ \text{for } n=0,1,2,\dots
    \label{radius_condition}
\end{align}
A special case of the above restriction, similar to condition \eqref{FUP_condition} (and also (1.3) in \cite{HelgeKnutsen_Daubechies_2020}), is
\begin{align}
    \pi R^2(n) \sim M^{\frac{n}{2}}.
    \label{radius_conditon_natural}
\end{align} 
Let $P_{n}(M,\mathcal{A})$ denote Daubechies' localization operator that projects onto the $n$-iterate spherically symmetric Cantor set with base $M$ and alphabet $\mathcal{A}$ defined in a disk whose radius $R(n)>0$ satisfies \eqref{radius_conditon_natural}. Then for some constant $B>0$, we find that 
\begin{align*}
    \| P_{n}(M,\mathcal{A})\|_{\text{op}}\leq B \left(\frac{|\mathcal{A}|}{M}\right)^{\frac{n}{2}} \ \ \text{for } n=0,1,2,\dots,
\end{align*}
where $|\mathcal{A}|<M$ denotes the size of the alphabet. This is similar to the asymptote for the mid-third Cantor set, which was estimated in part I \cite{HelgeKnutsen_Daubechies_2020} and was shown to be precise. Compared to the general case, we \textit{always} have that the asymptote is precise for the alphabet $\mathcal{A} =\overline{\mathcal{A}} =\{0,1,\dots,|\mathcal{A}|-1\}$ and \textit{never} precise for $\{M-1, M-2,\dots, M-|\mathcal{A}|\}$. Moreover,
for the alphabet $\overline{\mathcal{A}}$, the largest eigenvalue corresponds the Gaussian eigenfunction, independent of radius and iterate, which was not established for the mid-third Cantor set. 

While several properties transfer naturally from the non-indexed to the indexed case in section \ref{section_indexed_cantor_set}, for an arbitrary indexed base and alphabet, we do not always maintain good control over the operator norm. Therefore we present some simple conditions on the bases $\{M_j\}_j $ and alphabets $\{\mathcal{A}_j\}_j$ which guarantees that the operator norm decays to zero for increasing iterates $n$.

\section{Preliminaries}
\label{section_preliminaries}
\subsection{Daubechies' localization operator}
\label{Daubechies_localization_operator} 
In order to produce a \textit{joint} time-frequency representation of a signal $f\in L^2(\mathbb{R})$, we consider the \textit{Short-Time Fourier Transform} (STFT), dependent on a fixed, non-zero \textit{window} function $\phi: \mathbb{R}\to \mathbb{C}$. At point $(\omega, t)\in \mathbb{R}\times \mathbb{R}$, at frequency $\omega$ and time $t$, the STFT of $f$ with respect to the window $\phi$, is then given by 
\begin{align*}
    V_{\phi}f(\omega, t) = \int_{\mathbb{R}}f(x)\overline{\phi(x-t)}e^{-2\pi i \omega x}\mathrm{d}x.
\end{align*}
If we assume $\phi \equiv 1$, we retrieve the (regular) Fourier transform of $f$, without any time-dependence. For a joint time-frequency description, we assume $\phi$ to be non-constant. In particular, if we consider windows $\phi \in L^2(\mathbb{R})$, with $\|\phi \|_2 = 1$, the STFT becomes an isometry onto some subspace of $L^2(\mathbb{R}^2)$, i.e., $\langle V_{\phi} f, V_{\phi} g\rangle_{L^2(\mathbb{R}^2)} = \langle f, g\rangle$ for any $f,g\in L^2(\mathbb{R})$. Thus, we obtain a weakly defined inversion formula, where the original signal $f$ is recovered from $V_{\phi}f$ via inner products.

\textit{Daubechies' localization operator}, introduced in \cite{DaubechiesI.1988Tloa}, is based on the idea of modifying the STFT of $f$ by a multiplicative weight function $F(\omega,t)$ before recovering a time-dependent signal. The purpose of the weight function is to enhance and diminish different features of the $(\omega,t)$-domain $\mathbb{R}^2$, e.g., by projecting onto a subset of $\mathbb{R}^2$. Characterized by our choice of window $\phi$ and weight $F$, we denote the localization operator by $P_{F,\phi}$, which can be weakly defined as
\begin{align*}
    \langle P_{F,\phi} f, g \rangle := \langle F\cdot V_{\phi}f, V_{\phi} g \rangle_{L^2(\mathbb{R}^2)} \ \ \forall \ f,g\in L^2(\mathbb{R}). 
\end{align*}
Since its conception, this operator has not only been studied as an operator between $L^2$-spaces, but also more broadly as an operator between modulation spaces, and therein questions regarding boundedness and properties of eigenfunctions and egenvalues remain relevant (see \cite{Cordero_Grochenig_2003}, \cite{Cordero_Rodino_2008}, \cite{Bastianoni_Cordero_Nicola}). If we stick to the $L^2(\mathbb{R)}$-context, and assume the weight to be real-valued and integrable, the operator $P_{F,\phi}: L^2(\mathbb{R})\to L^2(\mathbb{R})$ becomes self-adjoint, compact. In particular, this means that the eigenfunctions of $P_{F,\phi}$ form an eigenbasis for $L^2(\mathbb{R})$, and the operator norm $\|P_{F,\phi}\|_{\text{op}}$ is given by the largest eigenvalue in absolute value. 

Similarly to Daubechies' classical paper \cite{DaubechiesI.1988Tloa} and what was done in part I \cite{HelgeKnutsen_Daubechies_2020}, we shall focus our attention to weights that are \textit{spherically symmetric}, that is, for some integrable function $\mathcal{F}:\mathbb{R}_{+}\to\mathbb{R}$, we consider
\begin{align}
    F(\omega,t) = \mathcal{F}(r^2), \ \text{where } r^2 = \omega^2 + t^2.
    \label{spherically_symmetric_weight}
\end{align}
Combined with a normalized Gaussian window, 
\begin{align}
    \phi(x) = 2^{1/4}e^{-\pi x^2},
    \label{normalized_Gaussian}
\end{align}
the eigenfunctions of $P_{F,\phi}$ are known, with explicit formulas for the associated eigenvalues:
\begin{theorem}
\textnormal{(Daubechies \cite{DaubechiesI.1988Tloa})} Let the weight $F$ and window $\phi$ be given by \eqref{spherically_symmetric_weight} and \eqref{normalized_Gaussian}, respectively. Then the eigenvalues of the localization operator $P_{F,\phi}$ read
\begin{align}
    &\lambda_k = \int_{0}^{\infty}\mathcal{F}\left(\frac{r}{\pi}\right)\frac{r^{k}}{k!}e^{-r}\mathrm{d}r, \ \text{for } k=0,1,2,\dots, 
    \label{thm_formula_eigenvalues_gen_spher_sym_weight}
\intertext{such that $P_{F,\phi}H_k = \lambda_k H_k$,
where $H_k$ denotes the $k$-th Hermite function,}
    &H_k(t) = \frac{2^{1/4}}{\sqrt{k!}}\left(-\frac{1}{2\sqrt{\pi}}\right)^{k}e^{\pi t^2}\frac{\mathrm{d}^k}{\mathrm{d}t^k}(e^{-2\pi t^2}).\nonumber
\end{align}
\label{theorem_derivation_spectrum_gaussian_window}
\end{theorem}

Interestingly, it was recently shown in \cite{Weyl_Heisenberg_Ensembles_Grochenig_2017} that any Hermite function $H_j$ as window and spherically symmetric weight yield localization operators with the same eigenbasis $\{H_k\}_k$. Another area that is being investigated is the aptly named \textit{inverse problem}, where one instead derives properties of the weight (symbol) based on knowledge of the eigenfuntions. E.g., for Daubechies' operators with a Gaussian window that project onto a simply connected domain $D\subseteq \mathbb{R}^2$, we know by \cite{Abreu_Dofler_Inverse_problem_2012} that if $H_j$ is an eigenfunction for some $j$, then $D$ reduces to a disk centered at the origin. More general situations are studied in \cite{Abreu_Grochenig_Romero_2016}, \cite{Abreu_Pereia_Romero_2017}.

Proceeding with the direct problem and Daubechies' classical result, formula \eqref{thm_formula_eigenvalues_gen_spher_sym_weight} represents a powerful tool for analyzing the localization operator and estimating the operator norm. In the subsequent discussion, we will only consider operators on the form as in Theorem \ref{theorem_derivation_spectrum_gaussian_window}. More precisely, we consider the case when $F$ projects onto some spherically symmetric subset $\pazocal{E}\subseteq \mathbb{R}^2$. That is, $\mathcal{F}(r) = \chi_{E}(r)$ for some subset $E\subseteq \mathbb{R}_{+}$ such that $\pazocal{E}= \{ (\omega, t)\in \mathbb{R}^2 \ | \ \omega^2+t^2\in E \}$. For the sake of simplicity, we will denote the associated localization operator by $P_{\pazocal{E}}$, whose eigenvalues are given by
\begin{align*}
    \lambda_{k}= \int_{\pi \cdot E}\frac{r^k}{k!}e^{-r}\mathrm{d}r \ \ \text{for } k=0,1,2,\dots, 
\end{align*}
where $\pi \cdot E := \{ x\in \mathbb{R}_{+} \ | \ x \pi^{-1} \in E\}$. Here it is worth noting that if we fix the measure $|E|$, we optimize the operator norm if $E$ corresponds to a ball centered at the origin. More precisely, 
\begin{align}
    \int_{\pi\cdot E}\frac{r^k}{k!}e^{-r}\mathrm{d}r \leq \int_{0}^{\pi|E|}e^{-r}\mathrm{d}r=1-e^{-\pi |E|} \ \ \text{for } k=0,1,2,\dots, 
    \label{optimal_localization_ball}
\end{align}
which was proved in Appendix A in \cite{HelgeKnutsen_Daubechies_2020}.

Since these will appear frequently, we shall denote the above integrands by
\begin{align*}
    f_k(r):= \frac{r^k}{k!}e^{-r} \ \text{for } k=0,1,2,\dots
\end{align*}
We recognize the function $f_k(r)$ as a gamma probability distribution, which is monotonically increasing for $r\in [0,k]$ and decreasing for $r\geq k$.

Finally, observe that these operators can also be studied from the perspective of Toeplitz operators on the Fock space. In particular, in \cite{Galbis_2021} Galbis considers Toeplitz operators with radial symbols and derive non-trivial norm estimates for such operators.

\subsection{Generalized Cantor set construction}
\label{section_cantor_set}
Fix a positive integer $M$, called the \textit{base}, and a non-empty proper subset $\mathcal{A}\subseteq \{0,1,\dots,M-1\}$, called the \textit{alphabet}.
The $n$-iterate (or $n$-order) \textit{discrete Cantor set} is then defined as a subset of $\mathbb{Z}_{M^n} = \{0,\dots,M^n-1\}$ of the form
\begin{align}
    \pazocal{C}_n^{(d)}(M,\mathcal{A}) := \left\{ \sum_{j=0}^{n-1} a_j M^j \ \big| \ a_j \in \mathcal{A} \ \text{for } j=0,1,\dots,n-1 \right\}.
    \label{discrete_Cantor_set}
\end{align}
The discrete version corresponds to the "\textit{continuous}" version based in the interval $[0,R]$ by
\begin{align}
    \pazocal{C}_n(R, M,\mathcal{A}) =  R M^{-n}\cdot \pazocal{C}_n^{(d)}(M, \mathcal{A})  + [0, R M^{-n}] \ \text{for } n=0,1,2,\dots,
    \label{continuous_Cantor_set_[0,R]}
\end{align}
where $a \cdot X = \{ x \ | \  x \cdot a^{-1}  \in X\}$ and $X + [0,b] = \cup_{x\in X} [x,x+b]$ for set $X$ and scalars $a,b >0$. Let $|\mathcal{A}|$ denote the cardinality of the alphabet $\mathcal{A}$. Then the measure of the $n$-iterate Cantor set is given by
\begin{align}
    |\pazocal{C}_n(R, M,\mathcal{A})| = \left(\frac{|\mathcal{A}|}{M}\right)^n R. 
\end{align}
In particular, for base $M = 3$ and alphabet $\mathcal{A} = \{0,2\}$, we recognize \eqref{continuous_Cantor_set_[0,R]} as the standard mid-third $n$-iterate Cantor set, with measure $(2/3)^n R$. Another noteworthy alphabet, that will appear frequently in the subsequent discussion, is
\begin{align*}
    \mathcal{A} = \overline{\mathcal{A}}:= \{0,1,\dots, |\mathcal{A}|-1\}, 
\end{align*}
which we will refer to as the \textit{canonical alphabet} of size $|\mathcal{A}|$. Such a redistribution of the alphabet does not alter the fractal dimension, i.e., the sets $\pazocal{C}_n(R,M,{\mathcal{A}})$ and $\pazocal{C}_n(R,M,\overline{\mathcal{A}})$  will still share the fractal dimension $\frac{\ln |\mathcal{A}|}{\ln M}$.  

For each $n$-iterate based in $[0,R]$, we define a corresponding map $\pazocal{G}_{n,R, M,\mathcal{A}}:\mathbb{R}\to [0,1]$, known as the \textit{Cantor function}, given by
\begin{align}
    \pazocal{G}_{n,R,M,\mathcal{A}}(x) := |\pazocal{C}_n(R,M,\mathcal{A})|^{-1}
    \begin{cases}
    0, x\leq 0,\\
    |\pazocal{C}_n(R,M,\mathcal{A})\cap[0,x]|, x>0.
    \end{cases}
    \label{Cantor_function}
\end{align}
Since $\pazocal{G}_{n,R,M,\mathcal{A}}(x) = \pazocal{G}_{n,1,M,\mathcal{A}}(xR^{-1})$, we set $\pazocal{G}_{n,M,\mathcal{A}} := \pazocal{G}_{n,1,M,\mathcal{A}}$, for simplicity. 
It is well-known that for the mid-third Cantor set the associated Cantor function is \textit{subadditive} (see \cite{JosefDobos_1996}). However, with an arbitrary alphabet, subadditivity can no longer be guaranteed. Instead we present a weaker version, sufficient for our purpose, utilizing the canonical alphabet (see Appendix \ref{section_weak_subadditivity_proof} for details).

\begin{lemma}
Let $\pazocal{G}_{n,M,\mathcal{A}}$ denote the Cantor-function, defined in \eqref{Cantor_function}. Then for any $x\leq y$,
\begin{align}
    \pazocal{G}_{n,M,\mathcal{A}}(y)-\pazocal{G}_{n,M,\mathcal{A}}(x) \leq \pazocal{G}_{n,M,\overline{\mathcal{A}}}(y-x).
    \label{weak_subadditivity_equation}
\end{align}
For $\mathcal{A} = \overline{\mathcal{A}}$, the above inequality is just standard subadditivity. 
\label{lemma_Cantor_function_weak_subadditivity_property}
\end{lemma}

For the disk of radius $R>0$, centered at the origin, we consider a spherically symmetric $n$-iterate, based on the $n$-iterate Cantor set in \eqref{continuous_Cantor_set_[0,R]}, as a subset of the form
\begin{align}
    \mathcal{C}_n(R, M,\mathcal{A}) = \{(\omega, t)\in \mathbb{R}^2 \ | \ \omega^2 + t^2 \in \pazocal{C}_n(R^2, M,\mathcal{A})\}\subseteq\mathbb{R}^2.
    \label{spherically_symmetric_Cantor_set}
\end{align}
This means we consider weights of the form 
\begin{align*}
    \mathcal{F}(r) = \chi_{\pazocal{C}_n(R^2, M,\mathcal{A})}(r) \ \text{for } R>0 \ \text{and } n=0,1,2,\dots
\end{align*}
The eigenvalues of Daubechies' localization operator $P_{\mathcal{C}_n(R, M,\mathcal{A})}$ then read
\begin{align}
    \lambda_k(\mathcal{C}_n(R, M,\mathcal{A})) = \int_{\pazocal{C}_n(\pi R^2, M,\mathcal{A})}f_k(r)\mathrm{d}r \ \text{for } k=0,1,2,\dots
    \label{eigenvalue_Cantor_set_formula}
\end{align}


\section{Localization on Generalized Spherically Symmetric Cantor set}
\label{localization_on_gen_spherically_symmetric_Cantor_set}
In this section we describe the behaviour of the operator norm, $\|P_{\mathcal{C}_n(R,M,\mathcal{A})}\|_{\text{op}}$ as a function of the iterates $n$. The results are formulated in section \ref{section_results_bounds}, with proofs and proof strategy in the subsequent sections \ref{section_main_tool_relative_areas}-\ref{section_counterexample_precise_asymptote}. 

\subsection{Results: Bounds for the operator norm}
\label{section_results_bounds}
Below we present three theorems regarding the operator norm of $P_{\mathcal{C}_n(R,M,\mathcal{A})}$.
The first theorem shows that $\|P_{\mathcal{C}_n(R,M,\mathcal{A})}\|_{\text{op}}$ can be bounded in terms of the "first eigenvalue" $\lambda_0(\mathcal{C}_n(R,M,\overline{\mathcal{A}}))$, thus revealing the significance of the canonical alphabet $\overline{\mathcal{A}} = \{0,1,\dots,|\mathcal{A}|-1\}$.
\begin{theorem}
\label{Theorem_operator_norm_eigenvalue_bound}
The operator norm of $P_{\mathcal{C}_n(R,M,\mathcal{A})}$ is bounded from above by 
\begin{align*}
    &\| P_{\mathcal{C}_n(R,M,\mathcal{A})}\|_{\text{op}}\leq 2 \lambda_{0}(\mathcal{C}_n(R,M,\overline{\mathcal{A}})).
\end{align*}
Further, for the canonical alphabet $\mathcal{A}=\overline{\mathcal{A}}$, we have that $\|P_{\mathcal{C}_n(R,M,\overline{\mathcal{A}})}\|_{\text{op}} = \lambda_{0}(\mathcal{C}_{n}(R,M,\overline{\mathcal{A}}))$.
\end{theorem}

In the next theorem we present an upper bound estimate for the operator norm $\| P_{\mathcal{C}_n(R,M,\mathcal{A})}\|_{\text{op}}$.  

\begin{theorem}
\label{Theorem_operator_norm_exponential_decay}
There exists a positive, finite constant $B$ only dependent on $|\mathcal{A}|$ and $M$ such that for each $n=0,1,2,\dots$
\begin{align*}
   \frac{\Big(\pi R^2 +1\Big)^{\frac{\ln|\mathcal{A}|}{\ln M}}}{|\mathcal{A}|^n\big(1-e^{-M^{-n}\pi R^{2}}\big)}\cdot \|P_{\mathcal{C}_n(R,M,\mathcal{A})}\|_{\text{op}} \leq B \ \ \ \forall \  \ \pi R^2 \in [0,M^n].
\end{align*}
\end{theorem}

Proofs of Theorem \ref{Theorem_operator_norm_eigenvalue_bound} and \ref{Theorem_operator_norm_exponential_decay} are found in section \ref{section_proof_thm1} and \ref{section_proof_thm2}, respectively.

\begin{remark}
If the alphabet, $\mathcal{A}$, is equal the canonical alphabet, $\overline{\mathcal{A}}$, then the left-hand-side of the inequality of Theorem \ref{Theorem_operator_norm_exponential_decay} can be bounded from below by a non-negative constant, thus making the asymptote \textit{precise}. 
\end{remark}

If we now enforce condition \eqref{radius_condition} on the radius $R$, we obtain the following result:

\begin{theorem}
\label{theorem_radius_conditon}
Suppose the radius $R$ depends on the iterates $n$ so that $R(n)\to\infty$ as $n\to \infty$ while $\pi R^2(n)\leq M^{n}$ for all $n=0,1,2,\dots$ Then there exist positive, finite constants $B_{L}\leq B_{U}$ only dependent on $M$ and $|\mathcal{A}|$ such that
\begin{enumerate}[label =  \textnormal{(\alph*)},itemsep=0.4ex, before={\everymath{\displaystyle}}]
\item for an arbitrary alphabet $\mathcal{A}$ we have the upper bound
\begin{align*}
    \| P_{\mathcal{C}_n(R(n), M, \mathcal{A})}\|_{\text{op}} \leq B_{U} \left(\frac{|\mathcal{A}|}{M}\right)^{n}\left(\pi R^2(n)\right)^{1-\frac{\ln|\mathcal{A}|}{\ln M}}, \ \ \text{and}
\end{align*}
\item for the canonical alphabet $\mathcal{A}=\overline{\mathcal{A}}$, we also have the lower bound
\begin{align*}
    \| P_{\mathcal{C}_n(R(n), M, \overline{\mathcal{A}})}\|_{\text{op}} \geq B_{L}\left(\frac{|\mathcal{A}|}{M}\right)^{n}\left(\pi R^2(n)\right)^{1-\frac{\ln|\mathcal{A}|}{\ln M}}. \ \ \ \ \ \ \ 
\end{align*}
\item Conversely, at any alphabet size $0<|\mathcal{A}|<M$, there exist alphabets $\mathcal{A}$ such that 
\begin{align*}
    \|P_{\mathcal{C}_n(R(n),M,\mathcal{A})}\|_{\text{op}}\cdot \|P_{\mathcal{C}_n(R(n),M,\overline{\mathcal{A}})}\|_{\text{op}}^{-1}\to 0 \ \ \text{as } n \to \infty.\ \
\end{align*}
\end{enumerate}
\end{theorem}

Recall that the quantity $\frac{\ln |\mathcal{A}|}{\ln M}$ is the fractal dimension of the Cantor set with base $M$ and alphabet $\mathcal{A}$. It should be noted that the exponent $\frac{\ln |\mathcal{A}|}{\ln M}-1$ in Theorem \ref{theorem_radius_conditon} (a) was already suggested in \cite{HelgeKnutsen_Daubechies_2020}. Since the associated upper bound holds for all alphabets and bases, this immediately begs the question whether the asymptote is precise \textit{regardless} of alphabet and base. By Theorem \ref{theorem_radius_conditon} (b) and also Corollary 4.1 in \cite{HelgeKnutsen_Daubechies_2020}, we conclude that the asymptote is precise for the Cantor set with canonical alphabet and for the mid-third Cantor set, respectively. However, by Theorem \ref{theorem_radius_conditon} (c), it becomes clear that we \textit{cannot} extend this result to every alphabet. A constructive proof of Theorem \ref{theorem_radius_conditon} (c) is found in section \ref{section_counterexample_precise_asymptote}.

\subsection{Main Tool: Relative Areas}
\label{section_main_tool_relative_areas}

Similarly to section 4 in \cite{HelgeKnutsen_Daubechies_2020}, our main tool is the concept of relative areas, namely
\begin{align}
    \pazocal{A}_{k,M,\mathcal{A}}(s,T) := \left[\sum_{a \in \mathcal{A}}\int_{s+aTM^{-1}}^{s+(a+1)TM^{-1}}f_k(r)\mathrm{d}r\right]\cdot\Bigg[\int_{s}^{s+T}f_k(r)\mathrm{d}r\Bigg]^{-1}.
    \label{relative_area_definition}
\end{align}
The relative areas measures the \textit{local} effect on the integrals that define the eigenvalues $\lambda_k(\dots)$ when we increase from one iterate $n$ to the next $n+1$. These are in general easier to work with rather than the eigenvalues themselves directly. Hence, we shall attempt to derive properties of the relative areas that transfer to the global behaviour of the eigenvalues. 

Initially, note that $\pazocal{A}_{0,M,\mathcal{A}}(s,T)$ is independent of the starting point $s\geq 0$, which yields the nice recursive relation for the first eigenvalue
\begin{align}
    \lambda_0(\mathcal{C}_{n+1}(R,M,\mathcal{A})) = \pazocal{A}_{0,M,\mathcal{A}}(\cdot, \pi R^2 M^{-n}) \lambda_0(\mathcal{C}_{n}(R,M,\mathcal{A})).
    \label{first_eigenvalue_alternative_formula}
\end{align}
For the canonical alphabet ${\mathcal{A}}=\overline{\mathcal{A}}$, \eqref{relative_area_definition} reduces to
\begin{align}
    \pazocal{A}_{k,M,\overline{\mathcal{A}}}(s,T) = \left[\int_{s}^{s+|\mathcal{A}|T M^{-1}}f_k(r)\mathrm{d}r\right]\cdot \Bigg[\int_{s}^{s+T}f_k(r)\mathrm{d}r\Bigg]^{-1}, 
    \label{relative_area_canonical_alphabet}
\end{align}
from which the relative area $\pazocal{A}_{0,M,\overline{\mathcal{A}}}(\cdot,T)$ attains the simple form
\begin{align}
     \pazocal{A}_{0,M,\overline{\mathcal{A}}}(\cdot,T) = \frac{1-e^{-\frac{|\mathcal{A}|}{M}T}}{1-e^{-T}}.
    \label{relative_area_simplified}
\end{align}
This relative area will play a significant role throughout the subsequent discussion. We conclude this section by showing $\mathcal{A}_{0,M,\overline{\mathcal{A}}}(\cdot,T)$ to be monotone, and illustrate how local effects can transfer to global behaviour.
\begin{lemma}
\label{lemma_eigenvalue_increasing_radius}
The relative area $\pazocal{A}_{0,M,\overline{\mathcal{A}}}(\cdot, T)$ is monotonically increasing in the argument $T$. 
\begin{proof}
Set $\epsilon := \frac{|\mathcal{A}|}{M}<1$, and simply differentiate \eqref{relative_area_simplified} to obtain
\begin{align*}
    \frac{\partial}{\partial T}\pazocal{A}_{0,M,\overline{\mathcal{A}}}(\cdot,T) = \frac{\partial}{\partial T}\left(\frac{1-e^{-\epsilon T}}{1-e^{-T}}\right) = e^{-T}\Big(\epsilon e^{(1-\epsilon)T}+(1-\epsilon)e^{-\epsilon T}-1\Big) \Big(1-e^{-T}\Big)^{-2}.
\end{align*}
It suffices to show that the \textit{second} factor in the above expression is always positive, i.e.,
\begin{align*}
    h_{\epsilon}(T) := \epsilon e^{(1-\epsilon)T}+(1-\epsilon)e^{-\epsilon T}-1\geq 0 \ \ \forall \ \ T\geq 0. 
\end{align*}
The latest claim is evident as 
\begin{align*}
    \frac{\partial h_{\epsilon}}{\partial T}(T) = \epsilon(1-\epsilon) e^{-\epsilon T}(e^{T}-1) \geq 0 \ \ \forall \ \ T\geq 0 \ \ \text{and } \ h_{\epsilon}(T) \to 0 \ \text{as } T\to 0. 
\end{align*}
\end{proof}
\end{lemma}

By monotonicity of $\pazocal{A}_{0,M,\mathcal{\overline{A}}}(\cdot, T)$ and the the recursive relation \eqref{first_eigenvalue_alternative_formula}, it follows that the associated eigenvalue $\lambda_{0}(\dots)$ is increasing as a function of the radius, i.e.,
\begin{align}
    \lambda_{0}(\mathcal{C}_{n}(R_1, M, \overline{\mathcal{A}}))\leq \lambda_{0}(\mathcal{C}_{n}(R_2, M, \overline{\mathcal{A}})) \ \ \forall \ R_1\leq R_2 \ \text{and } n=0,1,2,\dots
\end{align}
Monotonicity will also prove particularly useful both in section \ref{section_counterexample_precise_asymptote} and section \ref{section_indexed_cantor_set}.

\subsection{Proof of Theorem \ref{Theorem_operator_norm_eigenvalue_bound}}
\label{section_proof_thm1}

To begin with, we compare the relative areas with the canonical alphabet, for which we have the rather remarkable result. 

\begin{lemma}
\label{lemma_relative_area_canonical_alphabet_ordering}
Let $\{\pazocal{A}_{k,M,\overline{\mathcal{A}}}\}_k$ be given by \eqref{relative_area_canonical_alphabet}. Then for any $0<|\mathcal{A}|\leq M$ and $s\geq 0, T>0$, we have the ordering
\begin{align*}
    \pazocal{A}_{k,M,\overline{\mathcal{A}}}(s,T) \geq \pazocal{A}_{k+1,M,\overline{\mathcal{A}}}(s,T) \ \ \text{for } k=0,1,2,\dots
\end{align*}
\begin{proof}
For convenience, we define $\pazocal{A}_{k}(s,T,|\mathcal{A}|TM^{-1}) := \pazocal{A}_{k,M,\overline{\mathcal{A}}}(s,T)$, and show that the difference
\begin{align*}
    \pazocal{A}_{k}(s,T,t)-\pazocal{A}_{k+1}(s,T,t)\geq 0 \ \ \forall \ t\in [0,T].
\end{align*}
Firstly, we write the difference with a common denominator, which, by Fubini's theorem, yields
\begin{align*}
    &\int_{s}^{s+t}f_{k}(x)\mathrm{d}x \left[\int_{s}^{s+T}f_{k}(y)\mathrm{d}y\right]^{-1}-\int_{s}^{s+t}f_{k+1}(x)\mathrm{d}x \left[\int_{s}^{s+T}f_{k+1}(y)\mathrm{d}y\right]^{-1}\\
    &=\int_{s}^{s+t}\int_{s}^{s+T}\Big(f_{k}(x)f_{k+1}(y)-f_{k}(y)f_{k+1}(x)\Big)\mathrm{d}y\mathrm{d}x\cdot \left[\int_{s}^{s+T}f_k(r)\mathrm{d}r\int_{s}^{s+T}f_{k+1}(r)\mathrm{d}r\right]^{-1}.
\end{align*}
Since the integral
\begin{align*}
    \int_{s}^{s+t}\int_{s}^{s+t}\big(f_{k}(x)f_{k+1}(y)-f_{k}(y)f_{k+1}(x)\big)\mathrm{d}y\mathrm{d}x =0,
\end{align*}
the difference is reduced to
\begin{align*}
        \pazocal{A}_{k}(s,T,t)-\pazocal{A}_{k+1}(s,T,t)= \int_{s}^{s+t}\int_{s+t}^{s+T}\Big(f_{k}(x)f_{k+1}(y)-f_{k}(y)f_{k+1}(x)\Big)\mathrm{d}y\mathrm{d}x\cdot \Big[\dots\Big]^{-1}.
\end{align*}
Inserting the definition $f_k(r):= \frac{r^k}{k!}e^{-r}$ into the above integrand, we obtain
\begin{align*}
    f_{k}(x)f_{k+1}(y)-f_k(y)f_{k+1}(x) = \frac{1}{k!(k+1)!}x^k y^k(y-x)e^{-(x+y)}, 
\end{align*}
which is always positive since $s\leq x \leq s+t \leq y \leq s+T$.
\end{proof}
\end{lemma}

For the canonical alphabet, we can immediately conclude with the following corollary: 
\begin{corollary}
\label{corollary_canonical_alphabet_largest_eigenvalue}
The largest eigenvalue of the operator $P_{\mathcal{C}_{n}(R,M,\overline{\mathcal{A}})}$ is $\lambda_0(\mathcal{C}_{n}(R,M,\overline{\mathcal{A}}))$, and consequently the operator norm is given by
\begin{align*}
    \| P_{\mathcal{C}_{n}(R,M,\overline{\mathcal{A}})}\|_{\text{op}} = \lambda_0(\mathcal{C}_{n}(R,M,\overline{\mathcal{A}})).
\end{align*}
\begin{proof}
Since the relative area $\pazocal{A}_{0,M,\overline{\mathcal{A}}}(\cdot,T)$ is independent of the starting points $s\geq 0$ and bounds all $\{\pazocal{A}_{k,M,\overline{\mathcal{A}}}(s,T)\}_k$, it is clear that 
\begin{align*}
    &\lambda_k(\mathcal{C}_{n+1}(R,M,\overline{\mathcal{A}}))\leq \pazocal{A}_{0,M,\overline{\mathcal{A}}}(\cdot,\pi R^2 M^{-n})\ \lambda_k(\mathcal{C}_{n}(R,M,\overline{\mathcal{A}}))\\
    \intertext{which, by the relation \eqref{first_eigenvalue_alternative_formula} and observation \eqref{optimal_localization_ball}, eventually yields}
    &\lambda_k(\mathcal{C}_{n+1}(R,M,\overline{\mathcal{A}}))\leq \lambda_{k}(\mathcal{C}_{0}(R,M,\overline{\mathcal{A}}))\prod_{j=0}^{n}\pazocal{A}_{0,M,\overline{\mathcal{A}}}(\cdot, \pi R^2 M^{-j})\\
    & \ \ \ \ \ \ \ \ \ \ \ \ \ \ \ \ \ \ \ \ \ \ \ \ \ \leq \lambda_0(\mathcal{C}_{n+1}(R,M,\overline{\mathcal{A}})) \ \ \text{for } k=1,2,3\dots
\end{align*}
\end{proof}
\end{corollary}

For a general alphabet, we compare $\pazocal{A}_{0,M,\overline{\mathcal{A}}}(\cdot,T)$ to $\pazocal{A}_{k,M,\mathcal{A}}(s,T)$ for starting point $s\geq k$.

\begin{corollary}
\label{lemma_relative_area_s_geq_k}
Let $\{\pazocal{A}_{k,M,\mathcal{A}}\}_k$ be given by \eqref{relative_area_definition}. Then for $0<|\mathcal{A}|\leq M$ 
\begin{align*}
    \pazocal{A}_{0,M,\overline{\mathcal{A}}}(\cdot,T) \geq \pazocal{A}_{k,M,\mathcal{A}}(s,T) \ \forall \ s\geq k, \ T>0 \ \text{and } k=0,1,2,\dots
\end{align*}
\begin{proof}
Since $f_k(r)$ is monotonically decreasing for $r\geq k$, it follows that 
\begin{align*}
    \pazocal{A}_{k,M,\overline{\mathcal{A}}}(s,T) \geq \pazocal{A}_{k,M,\mathcal{A}}(s,T) \ \text{for } s\geq k,
\end{align*}
which combined with the ordering in Lemma \ref{lemma_relative_area_canonical_alphabet_ordering}, yields the result.
\end{proof}
\end{corollary}

By the same argument as in Corollary \ref{corollary_canonical_alphabet_largest_eigenvalue}, we utilize Corollary \ref{lemma_relative_area_s_geq_k} to obtain the bound 
\begin{align}
    \int_{\pazocal{C}_{n+1}(\pi R^2,M,\mathcal{A})+s}f_k(r)\mathrm{d}r \leq \pazocal{A}_{0,M,\overline{\mathcal{A}}}(\cdot,\pi R^2 M^{-n})\int_{\pazocal{C}_n(\pi R^2,M,\mathcal{A})+s}f_k(r)\mathrm{d}r&\nonumber\\ \leq \lambda_{0}(\mathcal{C}_{n+1}(R, M, \overline{\mathcal{A}})) \ \text{for } s\geq k&.
    \label{shifted_cantor_set_inequality}
\end{align}
Relating the shifted iterates $\pazocal{C}_n(\dots)+s$ to the non-shifted iterates $\pazocal{C}_n(\dots)$, we present an almost analogous statement to Lemma 3.5 in \cite{HelgeKnutsen_Daubechies_2020}.

\begin{lemma}
\label{lemma_non_shifted_iterates_weak_subadditivity}
Let $L>0$. Then for every fixed $k,n=0,1,2,\dots$, we have 
\begin{enumerate}[label =  \textnormal{(\Alph*)},itemsep=0.4ex, before={\everymath{\displaystyle}}]
    \item $\int_{\pazocal{C}_n(L,M,\mathcal{A})\cap[k,\infty[}f_k(r)\mathrm{d}r\leq \int_{\pazocal{C}_n(L,M,\overline{\mathcal{A}})+k}f_k(r)\mathrm{d}r\ $ and
    \item $\int_{\pazocal{C}_n(L,M,\mathcal{A})\cap [0,k]}f_k(r)\mathrm{d}r\leq \int_{\pazocal{C}_n(L,M,\overline{\mathcal{A}})+k}f_k(r)\mathrm{d}r$.
\end{enumerate}
\begin{proof}
Both cases (A) and (B) follow the same steps as in Lemma 3.5 in \cite{HelgeKnutsen_Daubechies_2020} but with the \textit{subadditivity} property for the mid-third Cantor function exchanged for the \textit{weaker version} of Lemma \ref{lemma_Cantor_function_weak_subadditivity_property} for the general Cantor function. Similarly to \cite{HelgeKnutsen_Daubechies_2020}, for the case (B), we use that $f(k+r)\geq f(k-r)$ for all $r\in [0,k]$ and reflect the elements $\pazocal{C}_n(L,M,\mathcal{A})\cap[0,k]$ about $r=k$ before proceeding with weak subadditivity. 

\end{proof}
\end{lemma}

Finally, by combining inequality \eqref{shifted_cantor_set_inequality} with Lemma \ref{lemma_non_shifted_iterates_weak_subadditivity}, we obtain the desired result 
\begin{align*}
    \lambda_k(\mathcal{C}_n(R,M,\mathcal{A})) \leq 2\lambda_0(\mathcal{C}_n(R,M,\overline{\mathcal{A}})) \ \ \text{for } k, n =0,1,2,\dots
\end{align*}

\subsection{Proof of Theorem \ref{Theorem_operator_norm_exponential_decay}}
\label{section_proof_thm2}

Proceeding, we present an \textit{exact} formula for the first eigenvalue: 

\begin{lemma}
\label{lemma_first_eigenvalue_formula}
The first eigenvalue of the operator $P_{\mathcal{C}_n(R,M,\mathcal{A})}$ is given by 
\begin{align}
    &\lambda_0(\mathcal{C}_n(R, M, \mathcal{A})) = \big(1-e^{- M^{-n}\pi R^2}\big)\prod_{j=1}^{n}\ \sum_{a_j\in \mathcal{A}}e^{-a_j M^{-j}\pi R^2}.
    \label{first_eigenvalue_equality}
\intertext{Further, we have the inequality}
    &\lambda_0(\mathcal{C}_n(R, M, \mathcal{A})) \leq \big(1-e^{- M^{-n}\pi R^2}\big)\prod_{j=1}^{n}\left(\frac{1-e^{-|\mathcal{A}|M^{-j}\pi R^2}}{1-e^{-M^{-j}\pi R^2}}\right),
    \label{first_eigenvalue_bound}
\end{align}
with equality precisely when $\mathcal{A} = \overline{\mathcal{A}} = \{0,1,\dots,|\mathcal{A}|-1\}$.
\begin{proof}
For simplicity, set $\pazocal{R} := \pi R^2$. By definition \eqref{relative_area_definition}, it is straightforward to compute the relative area $\pazocal{A}_{0,M,\mathcal{A}}$, which inserted into the recursive relation \eqref{first_eigenvalue_alternative_formula} yields
\begin{align*}
    &\lambda_0(\mathcal{C}_{n+1}(R,M,\mathcal{A}))
    = \left(\frac{1-e^{-\pazocal{R}M^{-(n+1)}}}{1-e^{-\pazocal{R}M^{-n}}}\right)\sum_{a_{n+1}\in\mathcal{A}}e^{-a_{n+1}\pazocal{R}M^{-(n+1)}} \lambda_0(\mathcal{C}_n(R,M,\mathcal{A})). 
\intertext{This in return means} 
    &\lambda_0(\mathcal{C}_{n+1}(R,M,\mathcal{A})) = \big(1-e^{-\pazocal{R}}\big) \prod_{j=1}^{n+1}\left(\frac{1-e^{-\pazocal{R}M^{-j}}}{1-e^{-\pazocal{R}M^{-(j-1})}}\right)\sum_{a_{j}\in\mathcal{A}}e^{-a_{j}\pazocal{R}M^{-j}}.
\end{align*}
Since the product $\prod_{j}\big(1-e^{-\pazocal{R}M^{-j}}\big)(1-e^{-\pazocal{R}M^{-(j-1)}}\big)^{-1}$ is telescoping, only the \textit{initial} denominator and \textit{final} numerator remain, and identity \eqref{first_eigenvalue_equality} readily follows. For the inequality case, merely note that the negative exponential function is monotonically decreasing and with the canonical alphabet $\overline{\mathcal{A}} = \{0,1,\dots,|\mathcal{A}|-1\}$, we recognize the appearance of the geometric series.   

\end{proof}
\end{lemma}


Using result \eqref{first_eigenvalue_bound} for the eigenvalue $\lambda_{0}(\mathcal{C}_n(R,M,\overline{\mathcal{A}}))$, we compute the asymptotes of said eigenvalue, which, combined with Theorem \ref{Theorem_operator_norm_eigenvalue_bound}, concludes the proof of Theorem \ref{Theorem_operator_norm_exponential_decay}.
\begin{proposition}
\label{proposition_first_eigenvalue_bounds}
There exists a finite constant $B\geq 1$ only dependent on $|\mathcal{A}|$ and $M$ such that
\begin{align*}
    B^{-1}\leq \frac{\Big(\pi R^2 +1\Big)^{\frac{\ln|\mathcal{A}|}{\ln M}}}{|\mathcal{A}|^n\big(1-e^{-M^{-n}\pi R^{2}}\big)}\cdot \lambda_{0}(\mathcal{C}_n(R,M,\overline{\mathcal{A}}))\leq B \ \ \ \forall \ \pi R^2 \in [0,M^n].
\end{align*}
\begin{proof}
The proof follows a similar technique to that of Proposition 4.1 in \cite{HelgeKnutsen_Daubechies_2020}, with minor modifications. By Lemma \ref{lemma_first_eigenvalue_formula}, the first eigenvalue of the operator $P_{\mathcal{C}_n(R,M,\overline{\mathcal{A}})}$, with the canonical alphabet $\overline{\mathcal{A}}$, is given by
\begin{align*}
    \lambda_0(\mathcal{C}_n(R, M, \overline{\mathcal{A}})) = \big(1-e^{- M^{-n}\pi R^2}\big)\prod_{j=1}^{n}\left(\frac{1-e^{-|\mathcal{A}|M^{-j}\pi R^2}}{1-e^{-M^{-j}\pi R^2}}\right).
\end{align*}
Hence, it suffices to show
\begin{align*}
    B^{-1} \leq \Big(x +1\Big)^{\frac{\ln|\mathcal{A}|}{\ln M}}|\mathcal{A}|^{-n}\prod_{j=1}^{n}\left(\frac{1-e^{-|\mathcal{A}|M^{-j}x}}{1-e^{-M^{-j}x}}\right) \leq B \ \ \forall \ x\in[0,M^n].
\end{align*}
Utilizing the factorization $(1-y^k) = (1-y)(1+y+y^2+\dots+y^{k-1})$ for $k\in\mathbb{N}$ and expressing the product as a sum, the above statement reads
\begin{align*}
    -\ln B \leq \sum_{j=1}^{n}\ln\left(1+e^{-x/M^j}+e^{-2x/M^j}+\dots + e^{-(|\mathcal{A}|-1)x/M^j}\right)&\\
    -\left(n-\frac{\ln(x+1)}{\ln M}\right)\ln|\mathcal{A}|\leq \ln B \ \ \ &\forall \ x\in[0,M^n]. 
\end{align*}
These inequalities follow by the aid of two claims
\begin{enumerate}[label =  (\roman*),itemsep=0.4ex, before={\everymath{\displaystyle}}]
    \item 
    \textit{there exists a finite, positive constant $\beta$ such that for $y\in[0,1]$}
    \begin{align*}
        -\beta\leq \sum_{j=1}^{\infty}\left[\ln\left(1+y^{1/M^j}+y^{2/M^j}+\dots +y^{(|\mathcal{A}|-1)/M^j}\right)-y^{1/M^j}\ln |\mathcal{A}|\right] \leq \beta, \ \textit{and}
    \end{align*}
    \item \textit{there exists a finite positive constant $\gamma$ such that}
    \begin{align*}
        -\gamma \leq \sum_{j=1}^{n}e^{-x/M^j}-\left(n-\frac{\ln(x+1)}{\ln M}\right) \leq \gamma \ \textit{for } x \in [0,M^n]. \ \ \ \ \ \
    \end{align*}
\end{enumerate}
For claim (i), we consider the function $\psi_{k}(y) := \ln(1+y+y^2+\dots+y^{k-1})-y \ln k$ for $y\in [0,1]$. Since $\psi_k$ is continuous and smooth in $[0,1]$, we have that
\begin{align*}
  c_k:=\max_{y\in [0,1]}\left|\frac{\mathrm{d}\psi_k(y)}{\mathrm{d}y}\right| < \infty,  
\end{align*}
from which we may define the linear spline
\begin{align*}
    h_k(y) := c_k\cdot\begin{cases}
    y,\ y \in [0,1/2]\\
    (1-y),\ y \in [1/2,1].
    \end{cases}
\end{align*} 
By the fact that $\psi_{k}(0) = \psi_{k}(1) =0$, it is clear that $\psi_k$ is bounded by the spline such that
\begin{align*}
    -h_k(y) \leq \psi_k(y) \leq h_k(y) \ \ \text{for } y\in[0,1].
\end{align*}
Thus, the sum in claim (i) is bounded from above and below by $\pm \sum_{j=0}^{\infty} h_{|\mathcal{A}|}(y^{1/M^j})$, respectively. From here the proof is essentially the same as for claim (i) in Proposition 4.1 in \cite{HelgeKnutsen_Daubechies_2020}. The same goes for claim (ii), where the only modification is that $3^j$ is exchanged for $M^j>1$.

\end{proof}
\end{proposition}

\subsection{Proof of Theorem \ref{theorem_radius_conditon} (c): Counterexample to precise asymptotic estimate}
\label{section_counterexample_precise_asymptote}
For our counterexample, we shall consider the reverse canonical alphabet, that is,
\begin{align*}
    \mathcal{A}=\underline{\mathcal{A}} := \{M-1,M-2,\dots,M-|\mathcal{A}|\},   
\end{align*}
where $M$ denotes the associated base of the Cantor set. As it turns out, this Cantor set construction is merely a shifted version of the Cantor set with canonical alphabet, $\overline{\mathcal{A}}$.   
\begin{lemma}
Let $\underline{\mathcal{A}}$ denote the reverse canonical alphabet of size $|\mathcal{A}|$. Then the Cantor set with alphabet $\underline{\mathcal{A}}$ and base $M$ can be expressed as
\begin{align*}
    \pazocal{C}_n(R,M,\underline{\mathcal{A}}) = R\cdot\big(M-|\mathcal{A}|\big)\sum_{j=1}^{n}M^{-j} +\pazocal{C}_n(R,M,\overline{\mathcal{A}}) \ \ \text{for } n=0,1,2,\dots
\end{align*}
\begin{proof}
By definition \eqref{discrete_Cantor_set}, the discrete Cantor set with the reverse canonical alphabet is given by
\begin{align*}
    \pazocal{C}_n^{(d)}(M,\underline{\mathcal{A}})&=\left\{ \sum_{j=0}^{n-1} a_j M^j \ \big| \ a_j = M-|\mathcal{A}|,\dots,M-1 \ \text{for } j=0,1,\dots,n-1 \right\}\\
    &=\left\{ \sum_{j=0}^{n-1} \big(a_j+M-|\mathcal{A}|\big) M^j \ \big| \ a_j = 0,1,\dots,|\mathcal{A}|-1 \ \text{for } j=0,1,\dots,n-1 \right\}\\
    &=\big(M-|\mathcal{A}|\big)\sum_{j=0}^{n-1}M^j + \pazocal{C}_{n}^{(d)}(M,\overline{\mathcal{A}}). 
\end{align*}
The continuous version then follows once we apply definition \eqref{continuous_Cantor_set_[0,R]} to the last identity.
\end{proof}
\end{lemma}

Hence, for the operator $P_{\mathcal{C}_{n}(R,M,\underline{\mathcal{A}})}$, the eigenvalues read 
\begin{align*}
    \lambda_k(\mathcal{C}_n(R,M,\underline{\mathcal{A}})) = \int_{\pi R_{\text{inner}}^2+\pazocal{C}_{n}(\pi R^2,M,\overline{\mathcal{A}})}f_k(r)\mathrm{d}r \ \ \text{for } k=0,1,2,\dots,
\end{align*}
where the inner radius $R_{\text{inner}} = R_{\text{inner}}(n,M,|\mathcal{A}|)\geq 0$ is given by
\begin{align}
       R_{\text{inner}}^2 = R^2\cdot \big(M-|\mathcal{A}|\big)\sum_{j=1}^{n}M^{-j}.
       \label{inner_radius}
\end{align}
According to condition \eqref{radius_condition}, we consider radii $R(n)$ that tends to infinity as $n\to\infty$, which also means that the inner radius $R_{\text{inner}} = R_{\text{inner}}(n)\to\infty$ as $n\to \infty$. Based on this inner radius, we may, in fact, exclude certain eigenvalues from being the largest.
More precisely, let $\floor{\cdot}$ denote the floor function, rounding down to the nearest integer. Since the difference between two integrands $f_{k+1}(r)-f_{k}(r)\geq 0$ for $r\geq k+1$, it is clear that largest eigenvalue $\lambda_k(\dots)$ must have index $k\geq \floor{\pi R_{\text{inner}}^2}$.

Proceeding, we consider the universal upper bound provided by Lemma \ref{lemma_non_shifted_iterates_weak_subadditivity}, namely
\begin{align}
     2\int_{k+\pazocal{C}_{n}(\pi R^2, M, \overline{\mathcal{A}})}f_k(r)\mathrm{d}r \geq \lambda_{k}(\mathcal{C}_n(R,M,\mathcal{A})) \ \ \text{for } k=0,1,2,\dots,
     \label{eigenvalue_universal_bound}
\end{align}
which holds for all alphabets, including $\mathcal{A}=\underline{\mathcal{A}}$.
We begin by computing the associated relative areas as $k\to\infty$.
\begin{lemma}
\label{lemma_limit_relative_area}
Let $\pazocal{A}_{k,M,\overline{\mathcal{A}}}(s,T)$ be the relative area given by \eqref{relative_area_canonical_alphabet} over the interval $[s,s+T]$. Suppose the starting point $s$ depends on $k$ so that $s=a k$ for $a>1$. Then  
\begin{align}
    \lim_{k\to \infty} \pazocal{A}_{k,M,\overline{\mathcal{A}}}(a k,T) = \frac{1-\exp\left[-\frac{|\mathcal{A}|}{M}T\left(1-\frac{1}{a}\right)\right]}{1-\exp\left[-T\left(1-\frac{1}{a}\right)\right]}.
    \label{limit_area}
\end{align}
For $s=k$, the limit reduces to $\lim_{k\to \infty}\pazocal{A}_{k,M,\overline{\mathcal{A}}}(k,T) = \frac{|\mathcal{A}|}{M}$.
\begin{proof}
Set $t:=|\mathcal{A}|TM^{-1}\in [0,T]$ for simplicity. By inserting the definition of the integrands $f_{k}(r) = \frac{r^k}{k!}e^{-r}$, we obtain 
\begin{align*}
    \lim_{k\to \infty} \pazocal{A}_{k,M,\overline{\mathcal{A}}}(a k,T) &= \lim_{k\to \infty}\int_{s}^{s+t}f_k(r)\mathrm{d}r\left[\int_{s}^{s+T}f_k(r)\mathrm{d}r\right]^{-1}\\
    &=\lim_{k\to \infty}\int_{0}^{t}e^{-r}\left(1+\frac{r}{ak}\right)^{k}\mathrm{d}r\left[\int_{0}^{T}e^{-r}\left(1+\frac{r}{ak}\right)^{k}\mathrm{d}r\right]^{-1}.
\end{align*}
Now, identity $\lim_{k\to \infty}\left(1+\frac{x}{k}\right)^{k} = e^{x}$ ensures we can utilize the dominated convergence theorem and exchange the order of the integrals and limit, so that
\begin{align*}
    \lim_{k\to \infty} \pazocal{A}_{k,M,\overline{\mathcal{A}}}(a k,T) &= \int_{0}^{t}e^{-r}\lim_{k\to \infty}\left(1+\frac{r}{ak}\right)^{k}\mathrm{d}r\left[\int_{0}^{T}e^{-r}\lim_{k\to \infty}\left(1+\frac{r}{ak}\right)^{k}\mathrm{d}r\right]^{-1}\\
    &=\int_{0}^{t}e^{-r\left(1-\frac{1}{a}\right)}\mathrm{d}r \left[\int_{0}^{T}e^{-r\left(1-\frac{1}{a}\right)}\mathrm{d}r\right]^{-1} = \frac{1-e^{-t\left(1-\frac{1}{a}\right)}}{1-e^{-T\left(1-\frac{1}{a}\right)}} \ \ \text{for }a>k.
\end{align*}
For $s=k$, i.e., $a=1$, the final integrands reduce to $e^{-r\left(1-\frac{1}{a}\right)}=1$.
\end{proof}
\end{lemma}

\begin{remark}
Comparing the above limit to the relative areas of $\lambda_0(\mathcal{C}_n(R,M,\overline{\mathcal{A}}))$ in \eqref{relative_area_simplified}, we find that $\lim_{k}\pazocal{A}_{k,M,\overline{\mathcal{A}}}(ak,T) = \pazocal{A}_{0, M,\overline{\mathcal{A}}}(\cdot, (1-\frac{1}{a})T)$.
\end{remark}
Although the limit of Lemma \ref{lemma_limit_relative_area} has a familiar form, it must be handled with some care as it, in fact, represents a lower bound rather than an upper one.
\begin{lemma}
\label{lemma_limit_area_caveat}
Suppose the factor $a\geq 1$. 
\begin{enumerate}[label =  \textnormal{(\Alph*)},itemsep=0.4ex, before={\everymath{\displaystyle}}]
    \item The limit relative area is an infimum in the sense that
\begin{align*}
    \pazocal{A}_{k,M,\overline{\mathcal{A}}}(ak,T)\geq \lim_{j\to \infty}\pazocal{A}_{j,M,\overline{\mathcal{A}}}(aj,T) \ \ \text{for } k=0,1,2,\dots \ \ \ \ \ \ \ \ \ \ \ \ 
\end{align*}
    \item Conversely, fix $\epsilon>0$. Then 
\begin{align*}
    \pazocal{A}_{k,M,\overline{\mathcal{A}}}(ak,T)\leq \lim_{j\to \infty}\pazocal{A}_{j,M,\overline{\mathcal{A}}}((1+\epsilon)j,T) \ \ \text{for } ak+T\leq (1+\epsilon)k.
\end{align*}
In particular, the inequality holds for $T\leq \frac{\epsilon k}{2}$ and $a\in\left[1,1+\frac{\epsilon }{2}\right]$.
\end{enumerate}
\begin{proof}
We start with (B) and relate this to (A). Set $t:=|\mathcal{A}|TM^{-1}$ and notice that
\begin{align*}
    \lim_{j\to\infty}\pazocal{A}_{j,M,\overline{\mathcal{A}}}((1+\epsilon)j,T)= \int_{ak}^{ak+t}e^{-x\left(1-\frac{1}{1+\epsilon}\right)}\mathrm{d}x\left[\int_{ak}^{ak+T}e^{-y\left(1-\frac{1}{1+\epsilon}\right)}\mathrm{d}y\right]^{-1} \ \ \textit{for } k=0,1,2,\dots
\end{align*}
On this form we compute the difference $\lim_j\pazocal{A}_{j,M,\overline{\mathcal{A}}}((1+\epsilon)j,T)-\pazocal{A}_{k,M,\overline{\mathcal{A}}}(ak,T)$, which, by the same technique as in Lemma \ref{lemma_relative_area_canonical_alphabet_ordering}, yields
\begin{align*}
    \int_{ak}^{ak+t}\int_{ak + t}^{ak + T}e^{-(x+y)}\left[y^k e^{\frac{x}{1+\epsilon}}-x^k e^{\frac{y}{1+\epsilon}}\right]\mathrm{d}y\mathrm{d}x \cdot\left[ \int_{ak}^{ak+T}y^{k}e^{-y}\mathrm{d}y\int_{ak}^{ak+T}e^{-\left(1-\frac{1}{1+\epsilon}\right)x}\mathrm{d}x\right]^{-1}.
\end{align*}
If the above integrand is \textit{always} positive or \textit{always} negative, this translates directly to the sign of the difference. This is the same as asking if the function 
\begin{align*}
    \ln Q(x,y) = \left(k\ln y-\frac{y}{1+\epsilon}\right)-\left(k\ln x +\frac{x}{1+\epsilon}\right) =: \ln q(y)-\ln q(x)
\end{align*}
is always positive or always negative for $ak \leq x\leq ak + t \leq y \leq ak +T$. Since the derivative $\left(\ln q\right)'(x)\geq 0$ for $x\in[k,(1+\epsilon)k]$, it follows that $\ln Q(x,y) \geq 0$ for $k\leq x\leq y\leq (1+\epsilon)k$ and thus part (B).  
For part (A), we simply consider $1+\epsilon = a$, which yields $(\ln q)'(x)\leq 0$ for $x\geq ak$.
\end{proof}
\end{lemma}

Now, it is tempting to suggest a common upper bound \textit{regardless} of starting point $s=ak$: 

\begin{example}
Suppose $L>1$ is constant and that $a \in [1,L]$ and $T\leq Lk$ for all iterates $n$. By Lemma \ref{lemma_limit_area_caveat} part (B), we have that
\begin{align*}
    \pazocal{A}_{k,M,\overline{\mathcal{A}}}(ak,T) \leq \lim_{j\to \infty}\pazocal{A}_{j,M,\overline{\mathcal{A}}}(2Lj,T)=\pazocal{A}_{0, M,\overline{\mathcal{A}}}\left(\cdot, \left(1-(2L)^{-1}\right)T\right).
\end{align*}
However, this estimate for the relative area is essentially the same as for $\lambda_{0}(\mathcal{C}_{n}(R,M,\overline{\mathcal{A}}))$, only with a scaled radius, $R \mapsto \sqrt{1-(2L)^{-1}} R$, which, by Theorem \ref{theorem_radius_conditon} (a) and (b), does not change the asymptotes.
\end{example}

With this example in mind, we divide the integral $\int_{k}^{\infty}f_k(r)\mathrm{d}r$ into a \textit{significant} and \textit{insignificant} part, thus, gaining more control over the starting points $s=ak$. In the next lemma, we show what we mean by insignificant. 

\begin{lemma}
\label{lemma_insignificant_integral}
For any fixed $\epsilon, \delta>0$, there exists a positive integer $K$, so that the integral
\begin{align*}
    \int_{(1+\epsilon)k}^{\infty}f_{k}(r)\mathrm{d}r = e^{-(1+\epsilon)k}\sum_{n=0}^{k}\frac{1}{n!}\big((1+\epsilon)k\big)^{n} < \delta \ \ \forall \ \ k\geq K.
\end{align*}
\begin{proof}
Using the lower bound version of \textit{Stirling's approximation formula} for the factorial
\begin{align*}
    \sqrt{2\pi} n^{n+\frac{1}{2}}e^{-n}\leq n! \ \ \text{for } n=1,2,3,\dots, 
\end{align*}
we determine a simplified upper bound for the summand 
\begin{align*}
    \sum_{n=1}^{k}\frac{1}{n!}\big((1+\epsilon)k\big)^{n} \leq \frac{1}{\sqrt{2\pi}}\sum_{n=1}^{k}\left(\frac{(1+\epsilon)k}{n}\right)^n \frac{e^{n}}{\sqrt{n}}\leq \frac{k}{\sqrt{2\pi}}(1+\epsilon)^{k}e^{k}.
\end{align*}
In the final inequality, we have used that the function 
\begin{align*}
    g_{A}(n) := \left(\frac{A}{n}\right)^{n} \ \ \text{is monotonically increasing for } 1\leq n \leq A e^{-1}.
\end{align*}
Applying the negative exponential $e^{-(1+\epsilon)k}$ to the upper bound, we obtain
\begin{align*}
    e^{-(1+\epsilon)k}\sum_{n=1}^{k}\frac{1}{n!}\big((1+\epsilon)k\big)^{n} \leq \frac{k}{\sqrt{2\pi}}\exp\big[-\big(\epsilon -\ln(1+\epsilon)\big)k\big], 
\end{align*}
and since $\epsilon-\ln(1+\epsilon)> 0$ for any fixed $\epsilon >0$, we are done. 
\end{proof}
\end{lemma}

\begin{proposition}
\label{proposition_significant_and_insignificant_integral_estimate}
Suppose the radius $R$ depends on the iterates $n$ such that $R(n) \to \infty$ as $n\to \infty$, and let $\gamma>0$ be a multiplicative constant. Then for every $\epsilon, \delta>0$, there exists a positive integer $N$ such that for every iterate $n\geq N$, we have that
\begin{enumerate}[label =  \textnormal{(\Alph*)},itemsep=0.4ex, before={\everymath{\displaystyle}}]
    \item $\int_{\big(\pazocal{C}_{n}(\pi R^2, M, \overline{\mathcal{A}})+k\big)\cap\big[k, \big(1+\frac{\epsilon}{2}\big)k\big]}f_{k}(r)\mathrm{d}r\leq 2\lambda_{0}(\mathcal{C}_n(\sqrt{\epsilon}R,M,\overline{\mathcal{A}})) \ \ \forall \ k\geq \gamma \cdot \pi R^2(n)$, and 
    \item $\int_{\big(\pazocal{C}_{n}(\pi R^2, M, \overline{\mathcal{A}})+k\big)\cap\big[ \big(1+\frac{\epsilon}{2}\big)k, \infty\big)}f_{k}(r)\mathrm{d}r \leq \delta\cdot \lambda_{0}(\mathcal{C}_n(R,M,\overline{\mathcal{A}})) \ \ \forall \ k\geq \gamma\cdot \pi R^2(n)$.
\end{enumerate}

\begin{proof}
We consider iterations $n$ so that the relevant indices $k\geq\gamma\cdot \pi R^2(n)\geq 1$.

For case (A), consider any fixed iterate $n_0$ so that
\begin{align*}
    M^{-n_0} \leq \frac{\gamma\cdot\epsilon}{2}.
\end{align*}
The interval size $|I|$ of the $n_0$-iterate Cantor set $\pazocal{C}_{n_0}(\pi R^2(n),M,\mathcal{A})$ then satisfies 
\begin{align*}
    |I| = \pi R^2(n) M^{-n_0} \leq \frac{\epsilon}{2} k \ \ \ \forall \ \ k\geq \gamma\cdot \pi R^2(n). 
\end{align*} 
By Lemma \ref{lemma_limit_area_caveat} (B) and identity $\lim_{k}\pazocal{A}_{k,M,\overline{\mathcal{A}}}(ak,T) =\pazocal{A}_{0,M,\overline{\mathcal{A}}}(\cdot,(1-a^{-1})T)$, we obtain
\begin{align}
    \int_{\big(\pazocal{C}_n(\dots)+k\big)\cap\left[k,\left(1+\frac{\epsilon}{2}\right)k\right]} f_{k}(r)\mathrm{d}r \leq \int_{k}^{\infty} f_{k}(r)\mathrm{d}r \cdot\prod_{j=n_0}^{n-1}\pazocal{A}_{0,M,\overline{\mathcal{A}}}\left(\cdot, \frac{\epsilon}{1+\epsilon}\pi R^2 M^{-j}\right).
    \label{asymptote_not_precise_case_1_id1}
\end{align}
By the recursive relation \eqref{first_eigenvalue_alternative_formula} for $\lambda_{0}(\mathcal{C}_n(\dots))$, it is clear that
\begin{align}
    \int_{0}^{\frac{\epsilon}{1+\epsilon}\pi R^2 M^{-n_0}}f_{0}(r)\mathrm{d}r \cdot\prod_{j=n_0}^{n-1}\pazocal{A}_{0,M,\overline{\mathcal{A}}}\left(\cdot, \frac{\epsilon}{1+\epsilon}\pi R^2 M^{-j}\right)\leq \lambda_{0}(\mathcal{C}_n(\sqrt{\epsilon}R,M,\overline{\mathcal{A}})),
    \label{asymptote_not_precise_case_1_id2}
\end{align}
where we have simplified the argument $\sqrt{\epsilon(1+\epsilon)^{-1}}R\leq \sqrt{\epsilon}R$ by monotonicity of $\lambda_{0}(\mathcal{C}_n(\dots))$.
Since the radius $R(n)\to\infty$ as $n\to \infty$, there exists a positive integer $n_1$ so that the integral
\begin{align}
    2\int_{0}^{\frac{\epsilon}{1+\epsilon}\pi R^2(n) M^{-n_0}}f_{0}(r)\mathrm{d}r\geq 1 \geq \int_{k}^{\infty}f_{k}(r)\mathrm{d}r \ \ \forall \ \ n\geq n_1. 
    \label{asymptote_not_precise_case_1_id3}
\end{align}
Combining the three inequalities \eqref{asymptote_not_precise_case_1_id1}-\eqref{asymptote_not_precise_case_1_id3}, then yields (A) with any $N\geq \max\{n_0, n_1\}$. 

For case (B), observe first, by the subadditivity of the Cantor function $\pazocal{G}_{n,M,\overline{\mathcal{A}}}$, that
\begin{align*}
    \int_{\big(\pazocal{C}_n(\dots)+k\big)\cap\left[\left(1+\frac{\epsilon}{2}\right)k,\infty\right)}f_{k}(r)\mathrm{d}r\leq \int_{\pazocal{C}_n(\dots)+\left(1+\frac{\epsilon}{2}\right)k}f_{k}(r)\mathrm{d}r.
\end{align*}
Further, by Corollary  \ref{lemma_relative_area_s_geq_k} and the recursive relation \eqref{first_eigenvalue_alternative_formula} for $\lambda_{0}(\mathcal{C}_n(\dots))$, we have
\begin{align*}
    \int_{\pazocal{C}_n(\dots)+\left(1+\frac{\epsilon}{2}\right)k}f_{k}(r)\mathrm{d}r&\leq \int_{\left(1+\frac{\epsilon}{2}\right)k}^{\infty}f_{k}(r)\mathrm{d}r\prod_{j=0}^{n-1}\pazocal{A}_{0,M,\overline{\mathcal{A}}}(\cdot,\pi R^2 M^{-j})\\
    &=\int_{\left(1+\frac{\epsilon}{2}\right)k}^{\infty}f_{k}(r)\mathrm{d}r\cdot \left[\int_{0}^{\pi R^2}f_0(r)\mathrm{d}r\right]^{-1}\cdot\lambda_{0}(\mathcal{C}_n(R,M,\overline{\mathcal{A}})). 
\end{align*}
We now apply Lemma \ref{lemma_insignificant_integral}, whereas for any $\delta>0$ there exists a positive integer $n_2$ so that
\begin{align*}
    \int_{\left(1+\frac{\epsilon}{2}\right)k}^{\infty}f_{k}(r)\mathrm{d}r<\frac{\delta}{2} \ \ \text{and } \int_{0}^{\pi R^2(n)}f_{0}(r)\mathrm{d}r\geq \frac{1}{2} \ \ \ \forall \ \ k\geq \gamma\cdot \pi R^2(n) \ \ \text{and } n\geq n_2,
\end{align*}
from which (B) follows with any $N\geq n_2$. For both cases, we chose $N\geq \max\{n_0,n_1,n_2\}$.

\end{proof}
\end{proposition}

\begin{proposition}
Let $\overline{\mathcal{A}}$ and $\underline{\mathcal{A}}$ denote the canonical and reverse canonical alphabet, respectively. Suppose the radius $R$ depends on the iterates $n$ such that $R(n)\to\infty$ as $n\to\infty$ while $\pi R^2(n)\leq M^{n}$. Then the quotient 
\begin{align*}
    \|P_{\mathcal{C}_n(R(n),M,\underline{\mathcal{A}})}\|_{\text{op}}\cdot \|P_{\mathcal{C}_n(R(n),M,\overline{\mathcal{A}})}\|_{\text{op}}^{-1}\to 0 \ \ \text{as } n\to \infty.
\end{align*}
\begin{proof}
Utilizing the universal upper bound \eqref{eigenvalue_universal_bound}, we have that
\begin{align*}
    \| P_{\mathcal{C}_n(R,M,\underline{\mathcal{A}})}\|_{\text{op}}\leq \sup_{k\geq \floor{\pi R^2_{\text{inner}}(n)}} 2 \int_{\pazocal{C}_n(\pi R^2,M,\overline{\mathcal{A}})+k}f_{k}(r)\mathrm{d}r \ \ \text{for } n=0,1,2,\dots
\end{align*}
Now fix some $\epsilon, \delta>0$. Since $R_{\text{inner}}(n)\to \infty$ as $n\to\infty$, we can apply Proposition \ref{proposition_significant_and_insignificant_integral_estimate} (A) and (B), and conclude that for some threshold iterate $N_0$
\begin{align*}
    \| P_{\mathcal{C}_n(R,M,\underline{\mathcal{A}})}\|_{\text{op}}\leq 4\lambda_{0}(\mathcal{C}_{n}(\sqrt{\epsilon}R,M,\overline{\mathcal{A}})+2\delta\cdot\lambda_{0}(\mathcal{C}_{n}(R,M,\overline{\mathcal{A}})) \ \ \forall \ n\geq N_0. 
\end{align*}
To see why this inequality readily yields the desired result, we insert the estimates of Theorem \ref{theorem_radius_conditon} (a) and (b). In particular, there exist a constant $B>0$ (independent of $\epsilon$) and threshold iterate $N_1$ so that
\begin{align*}
    &\lambda_{0}(\mathcal{C}_{n}(R(n), M, \overline{\mathcal{A}}))\geq B^{-1} \left(\frac{|\mathcal{A}|}{M}\right)^{n}\left(\pi R^2(n)\right)^{1-\frac{\ln|\mathcal{A}|}{\ln M}} \ \ \text{and } \\  &\lambda_{0}(\mathcal{C}_{n}(\sqrt{\epsilon}R(n),M,\overline{\mathcal{A}})) \leq  B \left(\frac{|\mathcal{A}|}{M}\right)^{n}\epsilon^{1-\frac{\ln |\mathcal{A}|}{\ln M}} \left(\pi R^2(n)\right)^{1-\frac{\ln|\mathcal{A}|}{\ln M}} \ \ \forall \ n\geq N_1.
\end{align*}
From these two estimates, and the fact that $\| P_{\mathcal{C}_{n}(R(n),M,\overline{\mathcal{A}})}\|_{\text{op}} = \lambda_{0}(\mathcal{C}_n(R,M,\overline{\mathcal{A}}))$, we obtain
\begin{align*}
    \| P_{\mathcal{C}_{n}(R(n),M,\underline{\mathcal{A}})}\|_{\text{op}}\cdot \| P_{\mathcal{C}_{n}(R(n),M,\overline{\mathcal{A}})}\|_{\text{op}}^{-1}\leq B^2\big( 4 \epsilon^{1-\frac{\ln |\mathcal{A}|}{\ln M}}+2\delta\big) \ \ \forall \ n \geq \max\{N_0,N_1\}.
\end{align*}
Since $\epsilon, \delta$ are arbitrary, and the constant $B$ does not depend on either, we are done. 
\end{proof}
\end{proposition}

\section{Further Generalizations: Indexed Cantor set}
\label{section_indexed_cantor_set}
So far in our discussion the base and alphabet have been \textit{fixed} throughout the iterations. 
Suppose we now instead let these two quantities \textit{vary}, to obtain an even more general Cantor set construction. In particular, if we index the bases $\{M_{j}\}_{j=1}^{\infty}$ and the alphabets $\{\mathcal{A}_j\}_{j=1}^{\infty}$ according to each iteration, we obtain the discrete construction
\begin{align}
    \pazocal{C}_{n}^{(d)}(\{M_j\}_{j}, \{\mathcal{A}_{j}\}_j):= \left\{\sum_{j=1}^{n}a_j \prod_{l=1}^{j-1} M_{l} \ \big| \ a_j \in\mathcal{A}_j \ \text{for } j=1,2,\dots,n\right\}, 
\end{align}
which in return yields the continuous version
\begin{align}
    \pazocal{C}_{n}(R, \{M_j\}_{j}, \{\mathcal{A}_{j}\}_j) = R \big(M_{1}& M_{2} \dots M_{n}\big)^{-1} \cdot \pazocal{C}_{n}^{(d)}(\{M_j\}_{j}, \{\mathcal{A}_j\}_{j})\\ +&\big[0,R\big(M_{1}M_{2}\dots M_{n}\big)^{-1}\big] \ \text{for } n=0,1,2,\dots\nonumber,
\end{align}
where an empty product is, by convention, defined as $1$. We shall refer to the above construction, $\pazocal{C}_n(R,\{M_j\}_j, \{\mathcal{A}_j\}_j)$, as an \textit{indexed Cantor set}. 

Similarly to \eqref{spherically_symmetric_Cantor_set}, we can also define an $n$-iterate spherically symmetric, indexed Cantor set, $\mathcal{C}_{n}(R,\{M_j\}_j,\{\mathcal{A}_j\}_j)$, and hence define the localization operator $P_{\mathcal{C}_{n}(R,\{M_j\}_j,\{\mathcal{A}_j\}_j)}$. With condition \eqref{radius_conditon_natural} in mind, a natural restriction on the radius seems to be
\begin{align}
    \pi R^2(n) \leq \gamma\cdot(M_1 M_2\dots M_n)^{1/2},
    \label{radius_indexed_cantor_set_condition}
\end{align}
for some finite constant $\gamma>0$,
and it is this restriction we utilize when formulating our localization results. The results are found in section \ref{section_results_decay_conditions}, with proofs in section \ref{section_indexed_cantor_set_proof_set_up} and \ref{section_indexed_cantor_set_proofs}.

\subsection{Results: Sufficient decay conditions}
\label{section_results_decay_conditions}

As should be expected, for an arbitrary indexed Cantor set, we cannot guarantee that the operator norm of the associated localization operator decays exponentially or even converges to zero for increasing iterates. Below we present two theorems that reflect this. In the first theorem we present some sufficient conditions for which the operator norm indeed decays exponentially. 
\begin{theorem}
Suppose we have an indexed Cantor set $\mathcal{C}_{n}(R,\{M_j\}_{j},\{\mathcal{A}\}_{j})$ such that
\begin{align*}
    \frac{|\mathcal{A}_j|}{M_j}\leq \epsilon < 1 \ \ \forall \ j\in\mathbb{N} \ \ \text{and } \ M_j\in [M,M^{1+\delta}]\ \text{for some finite } \delta>0 \ \text{and } M > 1. 
\end{align*}
Further, suppose that the radius $R$ satisfies the condition $\pi R^2(n) \leq \gamma\cdot \big(M_{1}M_{2}\dots M_{n}\big)^{1/2}$ for some finite $\gamma>0$. Then there exist finite constants $\alpha, \beta >0$ only dependent on $\epsilon, \delta$ and $\gamma$ such that the operator norm satisfies
\begin{align*}
    \| P_{\mathcal{C}_{n}(R(n), \{M_j\}_j, \{\mathcal{A}_{j}\}_j)}\|_{\text{op}}\leq \alpha e^{-\beta n} \ \ \text{for } n=0,1,2,\dots 
\end{align*}
\label{theorem_indexed_cantor_set_exponential_decay}
\end{theorem}
The second theorem shows that bounded quotients such that $\sup_{j}\frac{|\mathcal{A}_j|}{M_j}< 1$ is not itself a sufficient condition for the localization operator to converge to zero in the operator norm. 

\begin{theorem}
Fix $\gamma>0$, and suppose the radius $R$ satisfies $\pi R^2(n) = \gamma\cdot\big(M_1 M_2\dots M_n\big)^{1/2}$. Then there exist indexed bases $\{M_j\}_j$ and alphabets $\{\mathcal{A}_j\}_j$ with $\frac{|\mathcal{A}_j|}{M_j}\leq \epsilon < 1 \ \ \forall \ j \in \mathbb{N}$ such that
\begin{align*}
    \liminf_{n\to\infty}\| P_{\mathcal{C}_{n}(R(n), \{M_j\}_j, \{\mathcal{A}_{j}\}_j)}\|_{\text{op}} >0.
\end{align*}
\label{theorem_indexed_cantor_set_counterexample}
\end{theorem}

\subsection{Set-up and simple example}
\label{section_indexed_cantor_set_proof_set_up}
Based on our analysis of its simpler sister-operator $P_{\mathcal{C}_n(R,M,\mathcal{A})}$, we begin by deducing some analogous results for $P_{\mathcal{C}_n(R,\{M_j\}_j,\{\mathcal{A}_j\}_j)}$. 
In particular, identity \eqref{first_eigenvalue_alternative_formula} can be replaced by
\begin{align}
    \lambda_0(\mathcal{C}_{n+1}(R,\{M_j\}_j, \{\mathcal{A}_j\}_j)) = \pazocal{A}_{0,M_{n+1},\mathcal{A}_{n+1}}(\cdot,\pi R^2 (M_1 M_2 \dots M_n)^{-1})&
    \label{indexed_cantor_set_first_eigenvalue_recursive_relation}
    \\ \cdot\lambda_0(\mathcal{C}_n(R,\{M_j\}_j, \{\mathcal{A}_j\}_j))&.\nonumber
\end{align}
Further, as a matter of additional indexing in the proofs, it follows that an analogue of Theorem \ref{Theorem_operator_norm_eigenvalue_bound} must hold for the indexed localization operator, namely 
\begin{align*}
    \| P_{\mathcal{C}_n(R,\{M_j\}_j,\{\mathcal{A}_j\}_j)}\|_{\text{op}}\leq 2\lambda_{0}(\mathcal{C}_n(R,\{M_j\}_j, \{\overline{\mathcal{A}_j}\}_j)) \ \ \text{for } n=0,1,2,\dots,
\end{align*}
By the above inequality, it suffices to establish decay conditions for the first eigenvalue with \textit{canonical} alphabets. Therefore, by the recursive relation \eqref{indexed_cantor_set_first_eigenvalue_recursive_relation}, the relative areas $\{\pazocal{A}_{0,M_j,\overline{\mathcal{A}_j}}\}_{j}$ are of particular interest. 
Recall from result \eqref{relative_area_simplified} that these relative areas only depend on the base $M_j$ and alphabet $\overline{\mathcal{A}}_j$ via the quotient $\frac{|\mathcal{A}_j|}{M_j}$. Hence, for simplicity, we set 
\begin{align*}
    \pazocal{A}_{\frac{|\mathcal{A}_j|}{M_j}}(T):= \pazocal{A}_{0,M_j,\overline{\mathcal{A}_j}}(\cdot,T).
    \label{relative_area_simplified1}
\end{align*}
Recall, by Lemma \ref{lemma_eigenvalue_increasing_radius}, that the relative area $\pazocal{A}_{\frac{|\mathcal{A}_j|}{M}}(T)$ is a \textit{monotonically increasing} function for $T>0$. This again easily translates to the first eigenvalue with canonical alphabets being an increasing function as a function of the radius, i.e.,
\begin{align*}
    \lambda_0(\mathcal{C}_n(R_1,\{M_j\}_j,\{\overline{\mathcal{A}_j}\}_j)) \leq \lambda_0(\mathcal{C}_n(R_2,\{M_j\}_j,\{\overline{\mathcal{A}}_j\})) \ \ \forall \ \ R_1 \leq R_2 \ \text{and } n=0,1,2,\dots
\end{align*}
Thus, in the proof of Theorem \ref{theorem_indexed_cantor_set_exponential_decay} (and in Theorem \ref{theorem_indexed_cantor_set_counterexample}), we need only consider the equality case of condition \eqref{radius_indexed_cantor_set_condition}. From the general formula
\begin{align}
    &\lambda_{0}(\mathcal{C}_n(R(n),\dots)) = \left(1-e^{-\pi R^2(n)}\right) \prod_{j=1}^{n}\pazocal{A}_{\frac{|\mathcal{A}_{j}|}{M_{j}}}\left(\pi R^2(n)\cdot(M_1\dots M_{j-1})^{-1}\right),
    \intertext{we shall consider the eigenvalue on the form}
    &\lambda_{0}(\mathcal{C}_n(R(n),\dots)) = \left(1-e^{-\gamma\cdot(M_1\dots M_n)^{1/2}}\right)\prod_{j=1}^{n} \pazocal{A}_{\frac{|\mathcal{A}_{j}|}{M_{j}}}\left(\gamma\cdot\left[\frac{M_{j}\dots M_{n}}{M_1\dots M_{j-1}}\right]^{1/2}\right).
    \label{eigenvalue_indexed_cantor_set_simplified_identity}
\end{align}

We conclude this section with a simple example where the operator norm $\|P_{\mathcal{C}_{n}(\dots)}\|_{\text{op}}$ does \textit{not} decay to zero as the iterate $n\to \infty$.

\begin{example}
(Cantor set of positive measure) 
The measure of the $n$-iterate indexed Cantor set is given by
\begin{align*}
    |\pazocal{C}_n(R, \{M_j\}_j, \{\mathcal{A}_j\}_j)| = R \prod_{j=1}^{n}\frac{|\mathcal{A}_j|}{M_j}. 
\end{align*}
Hence, if $\frac{|\mathcal{A}_j|}{M_j}\to 1$ as $j\to\infty$ at a sufficient rate, the above product does not converge to zero as $n\to \infty$, and so the Cantor set itself has positive measure. That is, for such a choice of bases $\{M_j\}_j$ and alphabets $\{\mathcal{A}_j\}_j$, there exists some $0<\rho<1$ for which $|\pazocal{C}_n(R, \{M_j\}_j, \{\mathcal{A}_j\}_j)| \geq R\cdot \rho>0$ for all $n=0,1,2,\dots$ Using the integral formula for the first eigenvalue as a lower bound for the operator norm, it follows that
\begin{align*}
    \| P_{\mathcal{C}_n(R,\{M_j\}_j,\{\mathcal{A}_j\}_j)}\|_{\text{op}} \geq\int_{\pazocal{C}_n(\pi R^2,\dots)}e^{-r}\mathrm{d}r \geq e^{-\pi R^2}|\pazocal{C}_n(\pi R^2,\dots)|\geq e^{-\pi R^2} R\cdot \rho \ \ \forall \ n\in \mathbb{N}.
\end{align*}
This shows that for a Cantor set of positive measure $\liminf_{n}\|P_{\mathcal{C}_n(R,\dots)}\|_{\text{op}}>0$ when the radius $R$ is \textit{fixed}. 

If we instead let $R = R(n)$ increase according to the iterates $n$, we can choose canonical alphabets $\mathcal{A}_j = \overline{\mathcal{A}_j}$, which guarantees that the first eigenvalue $\lambda_{k}(R(n),\dots)$ is also increasing. Hence, $\liminf_{n}\|P_{\mathcal{C}_{n}(R(n),M_j,\overline{\mathcal{A}_j})}\|_{\text{op}}>0$.
\end{example}

\subsection{Proofs of Theorem \ref{theorem_indexed_cantor_set_exponential_decay} and Theorem \ref{theorem_indexed_cantor_set_counterexample}}
\label{section_indexed_cantor_set_proofs}
Since the first factor in \eqref{eigenvalue_indexed_cantor_set_simplified_identity} tends to $1$ as $n\to\infty$, convergence of the operator norm is determined by the remaining product $\prod_{j=1}^{n}(\dots)$. For simpler notation we assume the constant $\gamma=1$ in both proofs. 
\begin{proof}
(Theorem \ref{theorem_indexed_cantor_set_exponential_decay}) Initially, note that by the  interpretation as a relative area, the function $\pazocal{A}_{\frac{|\mathcal{A}_j|}{M_j}}$ must also be monotonically increasing in terms of the index $\frac{|\mathcal{A}_j|}{M_j}$. Since the index is bounded by $\epsilon$, we obtain the basic inequality
\begin{align*}
    \prod_{j=1}^{n}\pazocal{A}_{\frac{|\mathcal{A}_j|}{M_j}}\left(\left[\frac{M_{j}\dots M_{n}}{M_1\dots M_{j-1}}\right]^{1/2}\right)\leq \prod_{j=1}^{n}\pazocal{A}_{\epsilon}\left(\left[\frac{M_{j}\dots M_{n}}{M_1\dots M_{j-1}}\right]^{1/2}\right).
\end{align*}
From the second condition $M_j\in [M, M^{1+\delta}]$, we determine an upper bound for each argument in the relative area $\pazocal{A}_{\epsilon}(\dots)$. More precisely, we maximize the factors in the numerator, $M_{j},\dots, M_n \leq M^{1+\delta}$, and minimize the factors in the denominator, $M_{1},\dots, M_{j-1}\geq M$, so that
\begin{align*}
    \frac{M_{j}\dots M_{n}}{M_{1}\dots M_{j-1}} \leq M^{(1+\delta)(n-(j-1))}\cdot M^{-(j-1)} = M^{(1+\delta)n-(2+\delta)(j-1)} \ \ \text{for } j=1,2,\dots, n. 
\end{align*}
By monotonicity of $\pazocal{A}_{\epsilon}$, it follows that 
\begin{align*}
    \prod_{j=1}^{n}\pazocal{A}_{\frac{|\mathcal{A}_j|}{M_j}}\left(\left[\frac{M_{j}\dots M_{n}}{M_1\dots M_{j-1}}\right]^{1/2}\right)\leq \prod_{j=0}^{n-1}\pazocal{A}_{\epsilon}\left(M^{[(1+\delta)n-(2+\delta)j]/2}\right).
\end{align*}
We now define the set of factors which yields \textit{negative} exponents in the argument $M^{[\dots]/2}$, i.e.,
\begin{align}
    S_{n} := \big\{j\in \{0,1,\dots,n-1\} \ \big| \ (1+\delta)n-(2+\delta)j\leq 0\ \big\}.
    \label{factors_negative_exponents_def}
\end{align}
For these factors the argument satisfies $M^{[\dots]/2}\leq 1$, which, by monotonicity, means
\begin{align*}
    \prod_{j\in S_{n}}\pazocal{A}_{\epsilon}\left(M^{[\dots]/2}\right) \leq \prod_{j\in S_{n}}\pazocal{A}_{\epsilon}(1) = \pazocal{A}_{\epsilon}(1)^{|S_{n}|},
\end{align*}
where $|S_{n}|$ denotes the cardinality of $S_{n}$.
Since the relative areas $\pazocal{A}_{\epsilon}(T)\in (0,1) \ \ \forall \ \ T>0$, the remaining factors in $\{0,1,\dots,n-1\}\setminus S_{n}$ can be bounded by $1$. Hence, in order to prove exponential decay in the operator norm, it suffices to verify that $|S_n|\geq a \cdot n$ for some constant $a>0$. The latest claim is easily verified by a closer inspection of definition \eqref{factors_negative_exponents_def}, where $S_n$ can be expressed as
\begin{align*}
    S_{n} = \left\{ \left\lceil\left(\frac{1+\delta}{2+\delta}\right)n\right\rceil, \left\lceil\left(\frac{1+\delta}{2+\delta}\right)n\right\rceil+1,\dots, n-1\right\}.
\end{align*}
Thus, we conclude that $|S_n|\sim \left(1-\left(\frac{1+\delta}{2+\delta}\right)\right)n$.

\end{proof}

\begin{proof}
(Theorem \ref{theorem_indexed_cantor_set_counterexample}) We will construct a sequence of bases $\{M_j\}_j$ and alphabet cardinalities $\{|\mathcal{A}_j|\}_j$, where the product does not converge to zero as the iterate $n$ tends to infinity: 

For any fixed $\epsilon \in (0,1)$, chose an initial base $M>1$ and alphabet size $|\mathcal{A}|<M$ such that $\frac{|\mathcal{A}|}{M}\leq \epsilon$. Set $M_1 := M$ and $|\mathcal{A}_1| := |\mathcal{A}|$, and for $j>1$ set $M_j := M_{1}M_{2}\dots M_{j-1}$ and $|\mathcal{A}_j|:= \frac{|\mathcal{A}|}{M}M_j\in \mathbb{N}$ so that the fraction $\frac{|\mathcal{A}_j|}{M_j}$ remains constant and equal to $\frac{|\mathcal{A}|}{M}=:\theta$. With this construction, the product in \eqref{eigenvalue_indexed_cantor_set_simplified_identity} reads
\begin{align}
    &\prod_{j=1}^{n}\pazocal{A}_{\theta}\left(\left[\frac{M_{j}\dots M_{n}}{M_1\dots M_{j-1}}\right]^{1/2}\right)=\pazocal{A}_{\theta}\left([M_1\dots M_n]^{1/2}\right)\prod_{j=2}^{n}\pazocal{A}_{\theta}\left([M_{j+1}\dots M_{n}]^{1/2}\right)\label{product_proof_counterexample_construction}\\
    &= \pazocal{A}_{\theta}\left([M_1\dots M_n]^{1/2}\right)\big(\pazocal{A}_{\theta}\left([M_{3}\dots M_{n}]^{1/2}\right)\dots\pazocal{A}_{\theta}\left([M_{n-1}M_{n}]^{1/2}\right)\pazocal{A}_{\theta}\left(M_{n}^{1/2}\right)\pazocal{A}_{\theta}\left(1\right)\big). \nonumber
\end{align}
Proceeding, we make the following two basic observations
\begin{enumerate}[label =  (\roman*),itemsep=0.4ex, before={\everymath{\displaystyle}}]
    \setcounter{enumi}{1}
    \item $\pazocal{A}_{\theta}(T) \geq 1-e^{-\theta T}$ (which is still monotonically increasing) and
    \item $M_{j} \geq M \ \ \forall \ \ j\in \mathbb{N}$. 
\end{enumerate}
Combining these two observations with product \eqref{product_proof_counterexample_construction}, then yields the lower bound
\begin{align*}
    \prod_{j=2}^{n}\pazocal{A}_{\theta}\left([M_{j+1}\dots M_{n}]^{1/2}\right) \geq \prod_{m=2}^{n}\left(1-e^{-\theta N^{m}}\right) \geq \prod_{m=2}^{\infty}\left(1-e^{-\theta N^{m}}\right),
\end{align*}
where we have defined $N:=M^{1/2}$ for simplicity. It remains to show that the right-hand-side of the above inequality does not converge to zero for any fixed $N>1$ and $\theta \in (0,1)$. Exchanging the product for a sum, the statement is equivalent to showing that
\begin{align*}
    \sum_{m=2}^{\infty}\ln\left(1-e^{-\theta N^{m}}\right)> -\infty \ \ \text{for any fixed } N>1 \ \ \text{and } \theta \in (0,1).  
\end{align*}
Utilizing the lower bound $\ln(1-x) \geq -\left(\frac{x}{1-x}\right)$ for $0<x<1$, we obtain
\begin{align*}
    \sum_{m=2}^{\infty}\ln\left(1-e^{-\theta N^{m}}\right)\geq -\sum_{m=2}^{\infty}\frac{e^{-\theta N^m}}{1-e^{-\theta N^m}}\geq -\left(1-e^{-\theta}\right)^{-1}\sum_{m=2}^{\infty}e^{-\theta N^m},
\end{align*}
where the last sum converges by comparison with the geometric series. 

\end{proof}

\appendix
\section{Omitted Proof in section \ref{section_cantor_set}: Weak Subadditivy}
\label{section_weak_subadditivity_proof}
\begin{lemma}
Let $\pazocal{G}_{n,M,\mathcal{A}}$ denote the Cantor-function, defined in \eqref{Cantor_function}. Then for any $x\leq y$,
\begin{align}
    \pazocal{G}_{n,M,\mathcal{A}}(y)-\pazocal{G}_{n,M,\mathcal{A}}(x) \leq \pazocal{G}_{n,M,\overline{\mathcal{A}}}(y-x).
    \label{weak_subadditivity_equation_appendix}
\end{align}
\begin{proof}
Initially, note that the quantity $\pazocal{G}_{n,M,\mathcal{A}}(y)-\pazocal{G}_{n,M,\mathcal{A}}(x)$ measures the portion of the $n$-iterate Cantor set $\pazocal{C}_n$ contained in the interval $I = [x,y]$. Since the $n$-iterate is based in the interval $[0,1]$, we only need to consider the case when $I\subsetneq [0,1]$.

By definition, we have that the $n$-iterate Cantor set $\pazocal{C}_n$ effectively partitions $[0,1]$ into $M^{n}$ subintervals $\{I_k^{n}\}_{k=1}^{M^{n}}$, which we will refer to as the $n$-blocks. These blocks will either belong to $\pazocal{C}_n$ or have zero intersection, i.e., $|I_k^n \cap \pazocal{C}_n| = |I_k^n| (= M^{-n}) \ \text{or } =0.$
By the recursive construction, it follows that $I_k^{n-1}\cap \pazocal{C}_n$ either contains an entire alphabet of $n$-blocks or no $n$-blocks at all, i.e., $|I_k^{n-1}\cap \pazocal{C}_n| = |\mathcal{A}| M^{-n} \ \text{or } =0.$ In general, 
\begin{enumerate}[label =  (\roman*),itemsep=0.4ex, before={\everymath{\displaystyle}}]
    \item an $(n-j)$-block $I_k^{n-j}$ satisfies $|I_k^{n-j}\cap \pazocal{C}_n| = |\mathcal{A}|^j M^{-n} \ \text{or } =0.$
\end{enumerate}

Based on the distribution of the canonical alphabet, we deduce that 
\begin{enumerate}[label =  (\roman*),itemsep=0.4ex, before={\everymath{\displaystyle}}]
    \setcounter{enumi}{1}
    \item for any positive integer $m$ and positive parameter $a$,  \begin{align*}
        \pazocal{G}_{n,M,\overline{\mathcal{A}}}(m M^{j-n}+a) = \min\left\{1,\pazocal{G}_{n,M,\overline{\mathcal{A}}}(m M^{j-n})+\pazocal{G}_{n,M,\overline{\mathcal{A}}}(a)\right\}.
    \end{align*}
\end{enumerate}
Furthermore, utilizing result (i) and normalizing the Cantor function,
\begin{enumerate}[label =  (\roman*),itemsep=0.4ex, before={\everymath{\displaystyle}}]
    \setcounter{enumi}{2}
    \item for $m\leq M$ consecutive $(n-j)$-blocks, 
    \begin{align*}
        \pazocal{G}_{n,M,\overline{\mathcal{A}}}(m M^{j-n}) = \min\{m, |\mathcal{A}|\}|\mathcal{A}|^{j-n}. \ \ \ \ \ \ \ \ \ \ \ \ \ \ \ \ \ \ \ \ \ \ \ \ \ \ \ \ \ \ 
    \end{align*}
\end{enumerate}
Since the measure of any interval $I\subsetneq[0,1]$ can be uniquely expressed as
\begin{align}
   |I| = a + \sum_{j=1}^{n-1}m_j M^{j-n}
   \label{interval_length}
\end{align}
for integers $0\leq m_j < M$ and parameter $0\leq a \leq M^{1-n}$, we are now able to derive a more explicit formula for $\pazocal{G}_{n,M,\overline{\mathcal{A}}}$. In particular, by recursive use of (i) and finally inserting result (ii), we obtain
\begin{align}
    \pazocal{G}_{n,M,\overline{\mathcal{A}}}(|I|) &= \pazocal{G}_{n,M,\overline{\mathcal{A}}}\big(a + \sum_{j=1}^{n-1}m_j M^{j-n}\big) \nonumber\\
    &= \min\left\{1, \pazocal{G}_{n,M,\overline{\mathcal{A}}}\big(a) + \sum_{j=1}^{n-1}\pazocal{G}_{n,M,\overline{\mathcal{A}}}(m_j M^{j-n}\big)\right\}\nonumber\\
    &= \min\left\{1, \min\{aM^{n},|\mathcal{A}|\}|\mathcal{A}|^{-n} + \sum_{j=1}^{n-1}\min\{m_j,|\mathcal{A}|\}|\mathcal{A}|^{j-n}\right\}\label{cantor_function_canonical_alphabet_identity}.
\end{align}

Thus, it remains to show that the left-hand side of inequality \eqref{weak_subadditivity_equation_appendix} is bounded by \eqref{cantor_function_canonical_alphabet_identity}.
For this purpose, we shall utilize the following results:
\begin{enumerate}[label =  (\roman*),itemsep=0.4ex, before={\everymath{\displaystyle}}]
    \setcounter{enumi}{3}
    \item An interval of size $|I|\leq M^{1-n}$ satisfies
    \begin{align*}
        \pazocal{G}_{n,M,\mathcal{A}}(y)-\pazocal{G}_{n,M,\mathcal{A}}(x) \leq \min\{|I|M^{n},|\mathcal{A}|\}|\mathcal{A}|^{-n}.
    \end{align*}
    \item An interval of size $|I|\leq M^{j-n}$ for $j\in\mathbb{N}$ contains at most $|\mathcal{A}|^{j}$ $n$-blocks, i.e.,
    \begin{align*}
        \pazocal{G}_{n,M,\mathcal{A}}(y)-\pazocal{G}_{n,M,\mathcal{A}}(x) \leq |\mathcal{A}|^{j-n}. \ \ \ \ \ \ \ \ \ \ \ \ \ \ \ \ \ \ \ \ \ 
    \end{align*}
\end{enumerate}
In the case (iv), observe that $I$ will at most intersect \textit{two} $(n-1)$-blocks, say $I_{k}$ and $I_{k+1}$. If both intersections of $I_{k}$ and $I_{k+1}$ with $\pazocal{C}_n$ is non-zero, then the alphabet of $n$-blocks will have the same distribution in each $(n-1)$-block. Hence, by (i),
\begin{align*}
    |I\cap \pazocal{C}_n| = |I\cap I_k \cap \pazocal{C}_n|+|(I-M^{1-n})\cap I_k \cap \pazocal{C}_n| \leq |I_k \cap \pazocal{C}_n| = |\mathcal{A}|M^{-n},
\end{align*}
which yields (iv) after normalization. Result (v) follows by a similar argument. 

Based on observation (v), consider the case when $(m-1)M^{j-n} < |I| \leq m M^{j-n}$ for some integer $0<m\leq M$. Set $x_0:=x$ and  $x_m:=y$ and partition the interval $I=[x,y]=\cup_{l=0}^{m-1}[x_l,x_{l+1}]$ such that $x_{l}< x_{l+1}$ and $x_{l+1}-x_{l}\leq M^{j-n}$. Then 
\begin{enumerate}[label =  (\roman*),itemsep=0.4ex, before={\everymath{\displaystyle}}]
    \setcounter{enumi}{5}
    \item $\pazocal{G}_{n,M,\mathcal{A}}(y)-\pazocal{G}_{n,M,\mathcal{A}}(x) = \sum_{l=0}^{m-1}\big(\pazocal{G}_{n,M,\mathcal{A}}(x_{l+1})-\pazocal{G}_{n,M,\mathcal{A}}(x_{l})\big)\leq \min\{m, |\mathcal{A}|\} |\mathcal{A}|^{j-n}.$
\end{enumerate}
Finally, for an arbitrary interval $I=[x,y]\subseteq [0,1]$, we utilize the partitioning associated with \eqref{interval_length} along with results (iv)-(vi) to conclude that $\pazocal{G}_{n,M,\mathcal{A}}(y)-\pazocal{G}_{n,M,\mathcal{A}}(x)$ is indeed bounded from above by \eqref{cantor_function_canonical_alphabet_identity}.

\end{proof}
\end{lemma}
\addcontentsline{toc}{section}{References}
\bibliographystyle{IEEEtran}
\bibliography{Daubechies_time_frequency_localization_operator_on_Cantor_type_sets_part_II}

\end{document}